\documentclass[11pt]{article}
\usepackage[hmargin=2cm,vmargin=2.5cm]{geometry}
\usepackage[T1]{fontenc}
\usepackage{amsmath,amssymb,amsfonts,amsthm,mathrsfs}
\usepackage{chngcntr}
\usepackage[bookmarks=true]{hyperref}
\usepackage{color}
\usepackage[latin1]{inputenc}
\usepackage[shortlabels]{enumitem}

\numberwithin{equation}{section}

\newtheorem{THM}{Theorem}[section]
\newtheorem{PROP}[THM]{Proposition}

\newtheorem{LEM}[THM]{Lemma}
\newtheorem{DEF}[THM]{Definition}
\newtheorem{REM}[THM]{Remark}

\newtheorem{assumption}{Assumption}

\newtheorem{assumptionH}{Assumption}

\def \ms {\medskip}
\def \bs {\bigskip}
\def \expect {\mathbb{E}}
\def \ci {{\cal I}}
\def \p {\mathbb{P}}
\def \P{\mathbb{P}}
\def \ici {i\in \ci}

\def \sp {\hspace{0.25cm}}
\def \fr {\forall}
\def \qed {\blacksquare}

\DeclareMathOperator*{\esssup}{ess\,sup}

\makeatletter

\newcommand{\rnf}{\renewcommand{\thefootnote}{\arabic{footnote}}}
\newcommand{\thankyou}[2]{\stepcounter{footnote}\footnotetext[#1]{#2}}

\title{\vspace{-3.8cm}
\ \hspace{-2.70in}
{\small {}
}
\\
\vspace{3cm}
On a switching control problem with c\`adl\`ag costs}

\author{
\rnf Said Hamad\`ene\footnotemark[1]
\and
\rnf H\'ector Jasso-Fuentes\footnotemark[2]\  \  \footnotemark[3]
\and
\rnf Yamid A. Osorio-Agudelo\footnotemark[2]
}

\date{}

\begin{document}
\maketitle


\thankyou{1}{D\'epartement de Math\'ematiques, Equipe Statistique, et Processus, Universit\'e du Maine, Avenue
Olivier Messiaen, 72085 Le Mans, Cedex 9, France \href{mailto:hamadene@univ-lemans.fr}{hamadene@univ-lemans.fr}.}
\thankyou{2}{Departamento de Matem\'aticas. CINVESTAV-IPN. A. Postal 14-740, Ciudad de M\'exico, 07000, M\'exico. \href{mailto:\{h jasso, yaosorio\}@math.cinvestav.mx}{\{hjasso, yaosorio\}@math.cinvestav.mx}.}
\thankyou{3}{Corresponding author.}

\centerline{\today}

\begin{abstract} \noindent
This work addresses a switching control problem under which the cost associated with the changes of regimes is allowed to have discontinuities in time. Our main contribution is to show several characterizations of the optimal cost function as well as the existence of $\varepsilon$-optimal control policies. As a by-product, we also study the existence and uniqueness of solutions of a system of backward stochastic differential equations whose barriers (or obstacles) are discontinuous (in fact of c\`adl\`ag type) and depend itself on the unknown solution. At the last part of the paper, we study the case when  an underlying diffusion is part of the dynamic of the system. In this special case, the optimal payoff becomes a weak solution of the HJB system of PDEs with obstacles which is of quasi-variational type. This paper is somehow a continuation of the papers \cite{djehiche2009finite,hamadene2013viscosity} that consider continuous costs.

\end{abstract}
\noindent 


\section{Introduction}\label{sec1}
\setcounter{equation}{0}

%

Among the family of optimal control problems we can highlight those whose control is applied on the discontinuities to the dynamic. A special type of these problems is the so-called \emph{optimal multiple switching problems} consisting in configuring the state of system according to doing changes of  \emph{regimes} (a.k.a. \emph{configurations}) allowed for the controller. The times on  which these changes are triggered are also part of the control, so the controller needs to apply a sequence, say $(\tau_n,\xi_n)$ such that at time $\tau_n$, he/she changes the state from the regime $\xi_{n-1}$ to $\xi_n$, $n\geq 1$. The objective for him/her is to find an optimal sequence like the one above that maximizes a certain total payoff.

This class of problems has been studied in the literature by several authors. For instance, Carmona and  Ludkovski \cite{carmona2008pricing} study this kind of problems in order to find management optimal strategies with the purpose to release a power plant that converts natural gas into electricity and hence to sell this commodity in the market. Doucet and Ristic \cite{doucet2002recursive} apply the switching control theory to problems of target tracking that are commonly used in aerospace and electronic systems. Trigeorgis \cite{trigeorgis1993real,trigeorgis1996real} relates this type of problems to real option theory. Perhaps the most studied switching control problem is when only two-modes are considered.
Several authors have put attention on this type of problems (see e.g., Brekke and Oksendal \cite{brekke1991high,brekke1994optimal}, Hamad\`ene and Jeanblanc \cite{hamadene2007starting}, Duckworth and Zervos \cite{duckworth2000problem}, among others).

During the last decade, the switching control problem has been extensively studied by several authors including \cite{carmona2008pricing, chassagneux2011note, djehiche2009finite, hamadene2007starting,  HZ10, hutang,TY}, etc. (see also the references therein). 

However all the aforementioned papers consider the cases where the switching costs are continuous. To the best of our knowledge the case where the switching costs are discontinuous has not been considered yet. This is the main objective of this work. 

In this paper, as for the continuous switching costs, we show that the optimal payoffs are given by either a solution of a system of reflected BSDEs with obstacles depending on the solution or equivalently a system of processes expressed  through their corresponding Snell envelopes. This solution is discontinuous in time. On the other hand, while an optimal strategy may not exist, a nearly optimal strategy of switching always exists. Finally in the Markovian framework of randomness, the previous system provides a viscosity solution in weak sense of the Hamilton-Jacobi-Bellman system of PDEs associated with the switching problem. This paper is somehow the extension of the references Djehiche et. al. \cite{djehiche2009finite} and Hamad\`ene and Morlais \cite{hamadene2013viscosity} when the switching costs are of c\`adl\`ag type.

The rest of the paper is organized as follows: After this introductory part, in section \ref{sec:2} we introduce our switching problem and provide a verification theorem that is very common in control theory. In section \ref{sec3}, we present the existence and uniqueness of the solution for the system of RBSDEs with interconnected c\`adl\`ag obstacles and whose coefficients (drift and obstacles) depend on the unknown solution. Finally, in Section \ref{sec4} we show that under a Markovian framework i.e. the dynamic of the system is also governed by a underlying diffusion process, our unique solution, obtained in section \ref{sec3}, provides a weak viscosity solution for a system of a quasi-variational inequality with interconnected obstacles. $\qed$
\subsection{Notation and terminology}

Let $(\Omega,\mathcal{F},\mathbb{P})$ be a fixed probability space and  $B=\left(B_t\right)_{t\leq T}$ a $d$-dimensional Brownian motion with completed 
natural filtration $\left(\mathcal{F}_{t}:=\sigma\left\{B_s,s\leq t\right\}\right)_{t\leq T}$ thus it satisfies the usual conditions, i.e., it is right continuous and complete. Associated to $\mathbb{P}$, we denote by $\mathbb{E}$ its respective expectation.
\medskip

\noindent Next let us consider the following elements:
\begin{itemize}[leftmargin=*]
 \item $\left|\cdot\right|$ will denote the Euclidean norm in $\mathbb{R}^l$, for some appropriate $l\in\mathbb{N}$.
	\item Given $\theta\in[0,T]$, $L^{2}(\mathcal{F}_\theta)$ is the set of random variables $\xi$, $\mathcal{F}_\theta$-measurable and such that $\expect\big[\left|\xi\right|^{2}\big]<\infty$.
	\item $\mathcal{P}$ denotes the $\sigma$-algebra on $\left[0,T\right]\times\Omega$ of $(\mathcal{F}_{t})_{t\le T}$-progressively measurable sets.
	\item $\mathcal{H}^{2,l}$ denotes the set of $\mathcal{P}$-measurable processes $w=\left(w_t\right)_{t\leq T}$ with values in $\mathbb{R}^{l}$ such that $\left\|w\right\|_{\mathcal{H}^{2,l}}:=\expect\big[\int^{T}_{0}\left|w_s\right|^{2}ds\big]^{\frac{1}{2}}<\infty$. If $l=1$, then we will simply write $\mathcal{H}^{2,1}=\mathcal{H}^{2}$.
\item $\mathscr{S}^{2}$ stands for the set of $\mathcal{P}$-measurable, c\`adl\`ag, $\mathbb{R}$-valued processes $w=\left(w_t\right)_{t\leq T}$ such that
	$\left\|w\right\|_{\mathscr{S}^{2}}:=\{\expect\big[\sup_{t\leq T}\left|w_t\right|^{2}\big]\}^{\frac{1}{2}}<\infty$.
\item A random variable $\tau$ defined on $\Omega$ and valued in $\mathbb{R}_{+}\cup \left\{+\infty\right\}$ is called a stopping time with respect to the filtration $(\mathcal{F}_t)_{t\leq T}$, or simply an $\mathcal{F}_t$-stopping time, if for all $t\in\mathbb{R}_{+}$ $\left\{\omega|\tau(\omega)\leq t\right\}\in \mathcal{F}_{t}$.
\item $\mathcal{I}=\left\{1,\ldots,q\right\}$ denotes the set of indexes so-called set of configurations, 
while the notation $\mathcal{I}^{-i}$ means $\mathcal{I}-\left\{i\right\}$.
\item The notation $D^{2}_{xx}\phi$ and $D_x\phi$ denote the Hessian matrix and the gradient vector of the function $\phi$, respectively. 
\end{itemize}
\section{Model definition and preliminary results}\label{sec:2}

Consider the stochastic processes $\psi_i\in\mathcal{H}^{2}$, $i\in\mathcal{I}$, and $g_{ik}\in\mathscr{S}^{2}$, $i\in\mathcal{I}$ and $k\in\mathcal{I}^{-i}$, together with a sequence
\begin{equation}\label{stra}
\mathcal{S}=\left(\tau_n,\xi_n\right)_{n\geq0}
\end{equation}
of non-decreasing $\mathcal{F}$-stopping times $\tau_n$, and random variables $\xi_n$ which are $\mathcal{F}_{\tau_n}$-measurable  with values in $\mathcal{I}=\left\{1,\ldots,q\right\}$, such that $\tau_0=0$, $\xi_0=i$ for some $i\in\mathcal{I}$.

Together with these elements, define the functional $J$ as follows:
\begin{equation}\label{a01}
J(\mathcal{S},i)= \expect\Bigg[\sum^{\infty}_{n=0}\displaystyle\int^{\tau_{n+1}}_{\tau_n}\psi_{\xi_n}\left(s\right)ds-\displaystyle\sum^{\infty}_{n=1}g_{\xi_{n-1}\xi_{n}}(\tau_{n})\mathbf{1}_{\left[\tau_n<T\right]}\Bigg],
\end{equation}
with $\psi_{\xi_n}:=\psi_j$, when $\xi_n=j$; and the same reasoning applies to $g_{\xi_{n-1}\xi_n}$, i.e.,
$g_{\xi_{n-1}\xi_n}=g_{ik}$ if $\xi_{n-1}=i$ and $\xi_n=k$.
\begin{DEF}
A sequence $\mathcal{S}=\left(\tau_n,\xi_n\right)_{n\geq0}$ defined as in (\ref{stra}) is called a strategy or switching control policy for the controller. Furthermore, we say that a strategy $\mathcal{S}$ is \textit{admissible} if it satisfies the following condition:
\begin{equation}
\mathbb{P}\left[\tau_n<T,\hspace{0.1cm}\forall n\geq 0\right]=0\notag.
\end{equation}
For each $i=1,\ldots,q$, denote by $\mathcal{A}_{i}$ the set of admissible strategies with the property of $\tau_0=0$, $\xi_0=i$.
\end{DEF}
The processes $\psi_i$ and $g_{ik}$ are usually called the payoff rate per unit of the time and the switching cost, respectively. We will impose a condition to the processes $g_{ik}$, $i\in\mathcal{I}$, $k\in\mathcal{I}^{-i}$ that will be considered throughout this paper
\begin{assumption}\label{assA}
There exists a constant $\gamma>0$ such that the processes $g_{ik}\ge \gamma$ $\mathbb{P}$-a.s.
\end{assumption}
\medskip

A finite horizon switching control problem with $q$-modes and initial configuration $\xi_0=i$ for $i\in\mathcal{I}$, consists in finding an admissible sequence $S^{*}=\left(\tau^{*}_{n},\xi^{*}_{n}\right)_{n\geq0}$ $\in\mathcal{A}_i$ such that
\begin{equation}\label{pcg}
J(\mathcal{S}^{*},i)=\sup\limits_{\mathcal{S}\in\mathcal{A}_i}J(\mathcal{S},i)=:J^{*}(i),
\end{equation}
where $J$ is the functional defined in (\ref{a01}).

There is also a weaker formulation of what we understand for optimal strategy, namely, we say that $\mathcal{S}^{*}\in\mathcal{A}_i$ is $\varepsilon$-optimal strategy if for all $\varepsilon>0$, we have 
\begin{equation}
J(\mathcal{S}^{*},i)\geq J^{*}(i)-\varepsilon.\notag 
\end{equation}    
We first provide an existence result of $q$-interconnected processes, which will be useful later on. 
\begin{THM}\label{te}
Consider $q$ processes $\psi_i\in\mathcal{H}^{2}$, $i\in\mathcal{I}$ and $q(q-1)$ processes $g_{ik}\in\mathscr{S}^{2}$, $i\in\mathcal{I}$, $k\in\mathcal{I}^{-i}$. Then, under Assumption (\ref{assA}), there exist $q$ $\mathbb{R}-$valued c\`adl\`ag processes $\left(Y^i_{\cdot}:=(Y^{i}_{t}\right)_{t\leq T}, i=1,\ldots,q)\in(\mathscr{S}^{2})^{q}$ satisfying: $\forall i\in \ci$
\begin{equation}\label{a}\p-a.s., \,\,\forall t\leq T,\,\,
Y^{i}_{t}=\esssup_{\tau\geq t} \expect\bigg[\displaystyle\int^{\tau}_{t}\psi_i\left(s\right)ds+\max_{k\in\mathcal{I}^{-i}}\big(Y^{k}_{\tau}-g_{ik}(\tau)\big)\mathbf{1}_{\left[\tau<T\right]}\Big|\mathcal{F}_t\bigg].
\end{equation}
\end{THM}
\emph{Proof.} For $i\in\mathcal{I}$, and any $0\leq t\leq T$, use the sequence $(Y^{i,n}_{t})_{n\geq0}$ defined by:
\begin{equation}
Y^{i,0}_{t}=\expect\bigg[\displaystyle\int^{T}_{t}\psi_i\left(s\right)ds\Big|\mathcal{F}_t\bigg],\notag
\end{equation}
and for $n\geq 1$,
\begin{equation}\label{aa0}
Y^{i,n}_{t}=\esssup_{\tau\geq t}\expect\bigg[\displaystyle\int^{\tau}_{t}\psi_i\left(s\right)ds+\max_{k\in\mathcal{I}^{-i}}\big(Y^{k,n-1}_{\tau}-g_{ik}(\tau)\big)\mathbf{1}_{\left[\tau<T\right]}\Big|\mathcal{F}_t\bigg].
\end{equation}
First note that the process $(Y^{k,0}_{t})_{t\leq T}$ is continuous for all $k\in\mathcal{I}$. Next since the process $g_{ik}$ is c\`adl\`ag, $(Y^{i,1}_{t})_{t\leq T}$ is also a c\`adl\`ag process, and thus by an induction procedure we have that for all $n\geq 1$, $Y^{i,n}_{t}$ is c\`adl\`ag too.

Let us prove now that, for $i\in\mathcal{I}$, the sequence $(Y^{i,n}_{\cdot})_{n\geq 0}$ converges increasingly and pointwisely $\mathbb{P}$-a.s. for any $0\leq t\leq T$ and in the norm $\mathcal{H}^{2}$ to a c\`adl\`ag process $Y^{i}_{\cdot}$. To begin with, for any $n\geq 1$ let us define $\mathcal{A}^{i,n}_{t}=\left\{\mathcal{S}=(\tau_m,\xi_m)_{m\geq 0}:\xi_0=i, \tau_0= t \hspace{0.3cm} \textrm{and} \hspace{0.3cm} \tau_{n+1}=T\right\}$, and let us prove that for $N$ fixed, $Y^{i,N}_{\cdot}$ can be characterized by 
\begin{equation}\label{a11}
Y^{i,N}_{t}= \esssup_{\mathcal{S}\in \mathcal{A}^{i,N}_{t}}\expect\Bigg[\sum^{N}_{j=0}\displaystyle\int^{\tau_{j+1}}_{\tau_j}\psi_{\xi_j}\left(s\right)ds-\displaystyle\sum^{N-1}_{j=0}g_{\xi_j\xi_{j+1}}(\tau_{j+1})\mathbf{1}_{\left[\tau_{j+1}<T\right]}\bigg|\mathcal{F}_t\Bigg].
\end{equation}
Since the processes $g_{ik}$ for $i,k\in \mathcal{I}$ are c\`adl\`ag, it is not obvious to use the same procedure as given in Djehiche. et al. \cite{djehiche2009finite} Proposition 3-(ii). In contrast, we shall consider the sequence of $\varepsilon$-stopping times $(\tau^{\varepsilon}_{n})_{n\geq 0}$ given by follows: $\tau^{\varepsilon}_{0}:=t$,
\begin{equation}
\tau^{\varepsilon}_{1}:=\inf\Big\{s\geq t: Y^{i,N}_{s}\leq \max_{k\in\mathcal{I}^{-i}}\big(Y_s^{k,N-1}-g_{ik}\left(s\right)\big)+\frac{\varepsilon}{2}\Big\}\wedge T\notag
\end{equation}
and for $2\leq n\leq N$,
\begin{equation}
\begin{array}{ll}
\tau^{\varepsilon}_{n}:=\inf\Big\{s\geq\tau^{\varepsilon}_{n-1}: Y^{\hat{\xi}_{n-1},N-n+1}_{s}\leq \max\limits_{k\in\mathcal{I}^{-\hat{\xi}_{n-1}}}\big(Y_s^{k,N-n}-g_{\hat{\xi}_{n-1}k}(s)\big)+\frac{\varepsilon}{2^{n}}\Big\}\wedge T.& \\ &\\
\tau^{\varepsilon}_{N+1}:=T,
\end{array}\notag
\end{equation}
where
\begin{itemize}
  \item $\hat{\xi}_0:=i$, \hspace{0.2cm} $\hat{\xi}_1:=\arg\max\limits_{k\in\mathcal{I}^{-i}}\left\{Y^{k,N-1}_{\tau^{\varepsilon}_{1}}-g_{ik}(\tau^{\varepsilon}_{1})\right\}$
	\end{itemize}
	and for $n\geq 2$,
	\begin{itemize}

	\item $\hat{\xi}_n=\arg\max\limits_{k\in\mathcal{I}^{-\hat{\xi}_{n-1}}}\left\{Y^{k,N-n}_{\tau^{\varepsilon}_{n}}-g_{\hat{\xi}_{n-1}k}(\tau^{\varepsilon}_{n})\right\}.$
\end{itemize}
Note that by (\ref{aa0}) the process $(Y^{i,N}_{s}+\int^{s}_{t}\psi_i(r)dr)_{t\leq s\leq \tau^{\varepsilon}_{1}}$ is a super-martingale. Hence, if its Doob-Meyer decomposition is given by $(M_s-K_s)_{t\leq s\leq \tau^{\varepsilon}_{1}}$ (recall that $M$ is a martingale and $K$ a non-decreasing process), then by definition of $\tau^{\varepsilon}_{1}$, we have that $K_s=0$ for $s\in [t,\tau^{\varepsilon}_{1}]$, i.e., $(Y^{i,N}_{s}+\int^{s}_{t}\psi_i(r)dr)_{t\leq s\leq \tau^{\varepsilon}_{1}}$ is a martingale. Therefore,
\begin{equation}\label{a2}
\begin{array}{ll}
Y^{i,N}_{t} &=\expect\bigg[Y^{i,N}_{\tau^{\varepsilon}_{1}}+\displaystyle\int^{\tau^{\varepsilon}_{1}}_{t}\psi_i(r)dr\Big|\mathcal{F}_t\bigg] \\ \\ & \leq \expect\bigg[\max\limits_{k\in\mathcal{I}^{-i}}\left(Y_{\tau^{\varepsilon}_{1}}^{k,N-1}-g_{ik}\left(\tau^{\varepsilon}_{1}\right)\right)\mathbf{1}_{\left[\tau^{\varepsilon}_{1}<T\right]}+\frac{\varepsilon}{2}+\displaystyle\int^{\tau^{\varepsilon}_{1}}_{t}\psi_i(r)dr\Big|\mathcal{F}_t\bigg] \\ \\ & = \expect\bigg[\left(Y_{\tau^{\varepsilon}_{1}}^{\hat{\xi}_1,N-1}-g_{i\hat{\xi}_1}\left(\tau^{\varepsilon}_{1}\right)\right)\mathbf{1}_{\left[\tau^{\varepsilon}_{1}<T\right]}+\frac{\varepsilon}{2}+\displaystyle\int^{\tau^{\varepsilon}_{1}}_{t}\psi_i(r)dr\Big|\mathcal{F}_t\bigg] \\ \\ & =\expect\bigg[\displaystyle\int^{\tau^{\varepsilon}_{1}}_{t}\psi_i(r)dr-g_{i\hat{\xi}_1}\left(\tau^{\varepsilon}_{1}\right)\mathbf{1}_{\left[\tau^{\varepsilon}_{1}<T\right]}+Y_{\tau^{\varepsilon}_{1}}^{\hat{\xi}_1,N-1}\mathbf{1}_{\left[\tau^{\varepsilon}_{1}<T\right]}\Big|\mathcal{F}_t\bigg]+\frac{\varepsilon}{2}.
\end{array}
\end{equation}
Analogously, taking:
\begin{equation}
\tau^{\varepsilon}_{2}=\inf\Big\{s\geq\tau^{\varepsilon}_{1}, Y^{\hat{\xi}_1,N-1}_{s}\leq \max_{k\in\mathcal{I}^{-\hat{\xi}_1}}\big(Y_s^{k,N-2}-g_{\hat{\xi}_1k}(s)\big)+\frac{\varepsilon}{4}\Big\}\wedge T\notag
\end{equation}
we have again that $(Y^{\hat{\xi}_1,N-1}_{s}+\int^{s}_{\tau^{\varepsilon}_{1}}\psi_{\hat{\xi}_1}(r)dr)_{\tau^{\varepsilon}_{1}\leq s\leq \tau^{\varepsilon}_{2}}$ is a martingale. Arguing similarly as above, we have
\begin{equation}\label{a3}
\begin{array}{ll}
Y^{\hat{\xi}_1,N-1}_{\tau^{\varepsilon}_{1}}& =\expect\bigg[Y^{\hat{\xi}_1,N-1}_{\tau^{\varepsilon}_{2}}+\displaystyle\int^{\tau^{\varepsilon}_{2}}_{\tau^{\varepsilon}_{1}}\psi_{\hat{\xi}_1}(r)dr\Big|\mathcal{F}_{\tau^{\varepsilon}_{1}}\bigg]
  \\ \\ & \leq \expect\bigg[\max\limits_{k\in\mathcal{I}^{-\hat{\xi}_1}}\left(Y_{\tau^{\varepsilon}_{2}}^{k,N-2}-g_{\hat{\xi}_1k}\left(\tau^{\varepsilon}_{2}\right)\right)\mathbf{1}_{\left[\tau^{\varepsilon}_{2}<T\right]}+\frac{\varepsilon}{4}+\displaystyle\int^{\tau^{\varepsilon}_{2}}_{\tau^{\varepsilon}_{1}}\psi_{\hat{\xi}_1}(r)dr\Big|\mathcal{F}_{\tau^{\varepsilon}_{1}}\bigg] \\ \\ & =\expect\bigg[\left(Y_{\tau^{\varepsilon}_{2}}^{\hat{\xi}_2,N-2}-g_{\hat{\xi}_1\hat{\xi}_2}\left(\tau^{\varepsilon}_{2}\right)\right)\mathbf{1}_{\left[\tau^{\varepsilon}_{2}<T\right]}+\frac{\varepsilon}{4}+\displaystyle\int^{\tau^{\varepsilon}_{2}}_{\tau^{\varepsilon}_{1}}\psi_{\hat{\xi}_1}(r)dr\Big|\mathcal{F}_{\tau^{\varepsilon}_{1}}\bigg] \\ \\ & =\expect\bigg[\displaystyle\int^{\tau^{\varepsilon}_{2}}_{\tau^{\varepsilon}_{1}}\psi_{\hat{\xi}_1}(r)dr-g_{\hat{\xi}_1\hat{\xi}_2}\left(\tau^{\varepsilon}_{2}\right)\mathbf{1}_{\left[\tau^{\varepsilon}_{2}<T\right]}+Y_{\tau^{\varepsilon}_{2}}^{\hat{\xi}_2,N-2}\mathbf{1}_{\left[\tau^{\varepsilon}_{2}<T\right]}\Big|\mathcal{F}_{\tau^{\varepsilon}_{1}}\bigg]+\frac{\varepsilon}{4}.
\end{array}
\end{equation} 
Plugging (\ref{a3}) into (\ref{a2}), rearranging terms and since that $\left[\tau^{\varepsilon}_{2}<T\right]\subset \left[\tau^{\varepsilon}_{1}<T\right]$, we see that
\begin{equation}\label{a4}
Y^{i,N}_t\leq \expect\Bigg[\sum^{1}_{j=0}\bigg(\int^{\tau^{\varepsilon}_{j+1}}_{\tau^{\varepsilon}_j}\psi_{\hat{\xi}_j}(r)dr-g_{\hat{\xi}_j\hat{\xi}_{j+1}}\left(\tau^{\varepsilon}_{j+1}\right)\mathbf{1}_{\left[\tau^{\varepsilon}_{j+1}<T\right]}\bigg)+Y_{\tau^{\varepsilon}_{2}}^{\hat{\xi}_2,N-2}\mathbf{1}_{\left[\tau^{\varepsilon}_{2}<T\right]}\bigg|\mathcal{F}_t\Bigg]+\frac{\varepsilon}{2}+\frac{\varepsilon}{4}.\notag
\end{equation}
Repeating this procedure $N$ times, we obtain 
\begin{equation}\label{a5}
Y^{i,N}_t\leq\expect\Bigg[\sum^{N-1}_{j=0}\bigg(\int^{\tau^{\varepsilon}_{j+1}}_{\tau^{\varepsilon}_{j}}\psi_{\hat{\xi}_j}(s)ds-g_{\hat{\xi}_j\hat{\xi}_{j+1}}\left(\tau^{\varepsilon}_{j+1}\right)\mathbf{1}_{\left[\tau^{\varepsilon}_{j+1}<T\right]}\bigg)+Y_{\tau^{\varepsilon}_{N}}^{\hat{\xi}_N,0}\mathbf{1}_{\left[\tau^{\varepsilon}_{N}<T\right]}\bigg|\mathcal{F}_t\Bigg]+\varepsilon\left(\sum^{N}_{i=1}\frac{1}{2^{i}}\right).
\end{equation}
But
\begin{equation}\label{a6}
Y_{\tau^{\varepsilon}_{N}}^{\hat{\xi}_N,0}=\expect\bigg[\int^{T}_{\tau^{\varepsilon}_{N}}\psi_{\hat{\xi}_N}(s)ds\Big|\mathcal{F}_{\tau_N}\bigg].
\end{equation}
Plugging (\ref{a6}) into (\ref{a5}), and noting that $\left(\sum^{n}_{i=1}\frac{1}{2^{i}}\right)<1$, we deduce
\begin{equation}
Y^{i,N}_t\leq \expect\Bigg[\sum^{N}_{j=0}\int^{\tau^{\varepsilon}_{j+1}}_{\tau^{\varepsilon}_{j}}\psi_{\hat{\xi}_j}(s)ds-\sum^{N-1}_{j=0}g_{\hat{\xi}_j\hat{\xi}_{j+1}}\left(\tau^{\varepsilon}_{j+1}\right)\mathbf{1}_{\left[\tau^{\varepsilon}_{j+1}<T\right]}\bigg|\mathcal{F}_t\Bigg]+\varepsilon \hspace{0.4cm}\text{for all}\hspace{0.2cm}\varepsilon>0.\notag
\end{equation}
Since $(\tau^{\varepsilon}_{n},\hat{\xi}_n)_{0\leq n\leq N+1}$ belongs to $\mathcal{A}^{i,N}_{t}$, we can take essential supremum over $\mathcal{S}\in\mathcal{A}^{i,N}_{t}$ and then sending $\varepsilon\rightarrow 0$ to obtain 
\begin{equation}\label{a8}
Y^{i,N}_{t}\leq \esssup_{\mathcal{S}\in \mathcal{A}^{i,N}_{t}}\expect\Bigg[\sum^{N}_{j=0}\int^{\tau_{j+1}}_{\tau_{j}}\psi_{\hat{\xi}_j}(s)ds-\displaystyle\sum^{N-1}_{j=0}g_{\hat{\xi}_j\hat{\xi}_{j+1}}(\tau_{j+1})\mathbf{1}_{\left[\tau_{j+1}<T\right]}\bigg|\mathcal{F}_t\Bigg].
\end{equation}  
Now we derive the inverse inequality. Let $\mathcal{S}=(\tau_n,\xi_n)\in\mathcal{A}^{i,N}_{t}$ be an arbitrary strategy. Since $\tau_1\geq t$, $\mathbb{P}$-a.s., and $\xi_0=i$, then from (\ref{aa0}) we have 
\begin{equation}\label{a9}
\begin{array}{ll}
Y^{i,N}_{t} & =\esssup\limits_{\tau\geq t}\expect\bigg[\displaystyle\int^{\tau}_{t}\psi_i\left(s\right)ds+\max_{k\in\mathcal{I}^{-i}}\left(Y^{k,N-1}_{\tau}-g_{ik}(\tau)\right)\mathbf{1}_{\left[\tau<T\right]}\Big|\mathcal{F}_t\bigg] \\ \\ & \geq \expect\bigg[\displaystyle\int^{\tau_1}_{t}\psi_i\left(s\right)ds+\max_{k\in\mathcal{I}^{-i}}\left(Y^{k,N-1}_{\tau_1}-g_{ik}(\tau_1)\right)\mathbf{1}_{\left[\tau_1<T\right]}\Big|\mathcal{F}_t\bigg] \\ \\ & \geq \expect\bigg[\displaystyle\int^{\tau_1}_{t}\psi_i\left(s\right)ds+\left(Y^{\xi_1,N-1}_{\tau_1}-g_{i\xi_1}(\tau_1)\right)\mathbf{1}_{\left[\tau_1<T\right]}\Big|\mathcal{F}_t\bigg].
\end{array}
\end{equation}
In the same way, since $\tau_2\geq\tau_1$ and $\tau_1$ is also $\mathcal{F}_{\tau_2}$- measurable, then  
\begin{equation}
\begin{array}{ll}
Y^{\xi_1,N-1}_{\tau_1} & =\esssup\limits_{\tau\geq\tau_1}\expect\bigg[\displaystyle\int^{\tau}_{\tau_1}\psi_{\xi_1}\left(s\right)ds+\max_{k\in\mathcal{I}^{-\xi_1}}\left(Y^{k,N-2}_{\tau}-g_{\xi_1k}(\tau)\right)\mathbf{1}_{\left[\tau<T\right]}\Big|\mathcal{F}_{\tau_1}\bigg] \\ \\ & \geq \expect\bigg[\displaystyle\int^{\tau_2}_{\tau_1}\psi_{\xi_1}\left(s\right)ds+\max_{k\in\mathcal{I}^{-\xi_1}}\left(Y^{k,N-2}_{\tau_2}-g_{\xi_1k}(\tau_2)\right)\mathbf{1}_{\left[\tau_2<T\right]}\Big|\mathcal{F}_{\tau_1}\bigg]  \\ \\ & \geq \expect\bigg[\displaystyle\int^{\tau_2}_{\tau_1}\psi_{\xi_1}\left(s\right)ds+\left(Y^{\xi_2,N-2}_{\tau_2}-g_{\xi_1\xi_2}(\tau_2)\right)\mathbf{1}_{\left[\tau_2<T\right]}\Big|\mathcal{F}_{\tau_1}\bigg].
\end{array}\notag
\end{equation}
Plugging this last inequality into (\ref{a9}), rearranging terms and using that $\left[\tau_{2}<T\right]\subset \left[\tau_{1}<T\right]\in\mathcal{F}_{\tau_1}$, we see that
\begin{equation}
Y^{i,N}_{t}\geq\expect\Bigg[\displaystyle\sum^{1}_{j=0}\bigg(\int^{\tau_{j+1}}_{\tau_j}\psi_{\xi_j}\left(s\right)ds-g_{\xi_j\xi_{j+1}}(\tau_j)\mathbf{1}_{\left[\tau_{j+1}<T\right]}\bigg)+Y^{\xi_2,N-2}_{\tau_2}\mathbf{1}_{\left[\tau_2<T\right]}\bigg|\mathcal{F}_t\Bigg].\notag
\end{equation}
Continuing this procedure, we have
\begin{equation}
Y^{i,N}_{t}\geq\expect\Bigg[\displaystyle\sum^{N-1}_{j=0}\bigg(\int^{\tau_{j+1}}_{\tau_j}\psi_{\xi_j}\left(s\right)ds-g_{\xi_j\xi_{j+1}}(\tau_{j+1})\mathbf{1}_{\left[\tau_{j+1}<T\right]}\bigg)+Y^{\xi_N,0}_{\tau_N}\mathbf{1}_{\left[\tau_N<T\right]}\bigg|\mathcal{F}_t\Bigg].\notag
\end{equation} 
But again, since $Y^{\xi_N,0}_{\tau_N}=\expect\big[\int^{T}_{\tau_N}\psi_{\xi_n}(s)ds\big|\mathcal{F}_{\tau_N}\big]$, we get
\begin{equation}
Y^{i,N}_{t}\geq\expect\Bigg[\displaystyle\sum^{N}_{j=0}\int^{\tau_{j+1}}_{\tau_j}\psi_{\xi_j}\left(s\right)ds-\sum^{N-1}_{j=0}g_{\xi_j\xi_{j+1}}(\tau_{j+1})\mathbf{1}_{\left[\tau_{j+1}<T\right]}\bigg|\mathcal{F}_t\Bigg]\hspace{.3cm}\text{for all}\hspace{0.2cm}\mathcal{S}\in\mathcal{A}^{i,N}_{t}.\notag
\end{equation}
Thus, taking the essential supremum on $\mathcal{A}^{i,N}_{t}$, we get \begin{equation}
Y^{i,N}_{t}\geq\esssup_{\mathcal{S}\in \mathcal{A}^{i,N}_{t}}\expect\Bigg[\sum^{N}_{j=0}\int^{\tau_{j+1}}_{\tau_{j}}\psi_{\xi_j}(s)ds-\sum^{N-1}_{j=0}g_{\xi_j\xi_{j+1}}(\tau_{j+1})\mathbf{1}_{\left[\tau_{j+1}<T\right]}\bigg|\mathcal{F}_t\Bigg].\notag
\end{equation}
This last inequality together with (\ref{a8}), yield the characterization (\ref{a11}). 
\ms

\noindent Since $\mathcal{A}^{i,n}_{t}\subset\mathcal{A}^{i,n+1}_{t}$, we have $Y^{i,n}_{t}\leq Y^{i,n+1}_{t}$, $\mathbb{P}$-a.s. for all $t\in\left[0,T\right]$. On the other hand, by Assumption (\ref{assA}), we obtain for each $i\in\mathcal{I}$,
\begin{equation}\label{esti3.6}
 Y^{i,n}_{t}\leq \expect\bigg[\int^{T}_{t}\max_{[i=1,\ldots,q]}\left|\psi_i(s)\right|ds\Big|\mathcal{F}_t\bigg]\hspace{0.7cm} \text{for all}\hspace{0.2cm} t\leq T\hspace{0.2cm}\text{and}\hspace{0.2cm} n\geq 0 \notag
\end{equation}
and hence the sequence $(Y^{i,n}_{t})_{n \geq 1}$ is convergent. We now let $Y^{i}_{t}:=\lim_{n\rightarrow\infty}Y^{i,n}_{t}$ for $t\leq T$. Note that the process $Y^{i}_{\cdot}$ satisfies
\begin{equation}\label{a13}
Y^{i,0}_{t}\leq Y^{i}_{t}\leq \expect\bigg[\int^{T}_{t}\max_{[i=1,\ldots,q]}\left|\psi_i(s)\right|ds\Big|\mathcal{F}_t\bigg] \hspace{0.3cm}\text{for all}\hspace{0.2cm} t\leq T.
\end{equation}
Besides, $Y^{i}_{\cdot}$ is also c\`adl\`ag process. Indeed, from (\ref{aa0}) the process $(Y^{i,n}_{t}+\int^{t}_{0}\psi_i(s)ds)_{0\leq t\leq T}$ is a c\`adl\`ag super-martingale for all $i\in\mathcal{I}$ and $n\geq 1$. Thus its limit process $(Y^{i}_{t}+\int^{t}_{0}\psi_i(s)ds)_{0\leq t\leq T}$ is c\`adl\`ag as a limit of increasing sequence of c\`adl\`ag super-martingales (see Dellacherie and Meyer [\cite{dellacherie1980probabilites}, p. 86]), which gives the desired c\`adl\`ag property of $Y^{i}_{.}$. Moreover, from (\ref{a13}), the $L^{2}$-properties of $\psi_i$ and by Doob's maximal inequality, for each $i\in\mathcal{I}$, we have
\begin{equation}
\expect\bigg[\sup_{0\leq t\leq T}\left|Y^{i}_{t}\right|^{2}\bigg]<\infty\notag
\end{equation}
and hence by the Lebesgue dominated convergence theorem the sequence $(Y^{i,n}_{\cdot})_{n\geq 0}$ converges to $Y^{i}_{\cdot}$ in $\mathcal{H}^{2}$. Thus, by Snell envelope properties (see Proposition 2-(iv) in Djehiche, et al. \cite{djehiche2009finite}), the c\`adl\`ag processes $Y^{1}_{\cdot},\ldots,Y^{q}_{\cdot}$ satisfy (\ref{a}) since they are limits of the increasing sequence of processes $Y^{i,n}_{\cdot}$, for $i\in \mathcal{I}$, satisfying (\ref{a11}).$\hspace{1cm}\blacksquare$ 
\medskip

Let us show now some properties of the $\varepsilon$-strategy introduced in Theorem \ref{te}. 
\begin{PROP}\label{admissible}
The $\varepsilon$-strategy $\mathcal{S}^{\varepsilon}=\left(\tau^{\varepsilon}_n,\xi^{\varepsilon}_n\right)_{n\geq 0}$ defined as follows:
\begin{itemize}
	\item $\tau^{\varepsilon}_{0}:=0$, \hspace{.2cm} $\tau^{\varepsilon}_{1}:=\inf\Big\{s\geq 0: Y^{i}_{s}\leq \max\limits_{k\in\mathcal{I}^{-i}}\left(Y_s^{k}-g_{ik}\left(s\right)\right)+\frac{\varepsilon}{2}\Big\}\wedge T$
\end{itemize}
and, for $n\geq 2$,
\begin{itemize}
	\item $\tau^{\varepsilon}_{n}:=\inf\Big\{s\geq\tau^{\varepsilon}_{n-1}: Y^{\xi^{\varepsilon}_{n-1}}_{s}\leq \max\limits_{k\in\mathcal{I}^{-\xi^{\varepsilon}_{n-1}}}\left(Y_s^{k}-g_{\xi^{\varepsilon}_{n-1}k}(s)\right)+\frac{\varepsilon}{2^{n}}\Big\}\wedge T$
\end{itemize}
and the sequence $(\xi^{\varepsilon}_n)$ given by
\begin{itemize}
	\item $\xi^{\varepsilon}_{0}:=i$, \hspace{.2cm} $\xi^{\varepsilon}_1=\arg\max\limits_{k\in\mathcal{I}^{-i}}\left\{Y^{k}_{\tau^{\varepsilon}_{1}}-g_{ik}(\tau^{\varepsilon}_{1}) \right\}$
	\end{itemize}
	and for $n\geq 2$,
	\begin{itemize}
\item $\xi^{\varepsilon}_n=\arg\max\limits_{k\in\mathcal{I}^{-\xi_{n-1}}}\left\{Y^{k}_{\tau^{\varepsilon}_{n}}-g_{\hat{\xi}_{n-1}k}(\tau^{\varepsilon}_{n})\right\},$
\end{itemize}
is admissible.
\end{PROP}
\textit{Proof.} Suppose for contradiction that $\mathcal{S}^{\varepsilon}$ is not admissible, that is, $\P[\tau^{\varepsilon}_{n}<T,\text{for all}\hspace{0.2cm} n\geq 1]>0$. Then, by definition of $\tau^{\varepsilon}_{n}$ we have
\begin{equation}
\P\left[Y^{\xi^{\varepsilon}_{n-1}}_{\tau^{\varepsilon}_{n}}\leq Y^{\xi^{\varepsilon}_{n}}_{\tau^{\varepsilon}_{n}}-g_{\xi^{\varepsilon}_{n-1}\xi^{\varepsilon}_{n}}(\tau^{\varepsilon}_{n})+\frac{\varepsilon}{2^{n}},\hspace{0.2cm} \xi^{\varepsilon}_{n}\in\mathcal{I}^{-\xi^{\varepsilon}_{n-1}}, \forall n\geq 1\right]>0.\notag
\end{equation}
If the event $B=\left\{\omega\in\Omega:\tau^{\varepsilon}_{n}(\omega)<T, \forall n\geq 1\right\}$ has positive probability, then  there is a state $i_1\in\mathcal{I}$ and a loop $i_1,i_2\ldots,i_k$ (with $i_1=i_k$) of elements of $\mathcal{I}$ (recall that $\mathcal{I}$ is a finite set), and subsequence $\tau^{\varepsilon}_n,\ldots,\tau^{\varepsilon}_{n+k}$ corresponding of this configuration such that  
\begin{equation}\label{limite}
\P\left[Y^{i_l}_{\tau^{\varepsilon}_{n+l}}\leq Y^{i_{l+1}}_{\tau^{\varepsilon}_{n+l}}-g_{i_l,i_{l+1}}(\tau^{\varepsilon}_{n+l})+\frac{\varepsilon}{2^{n}}, l=1,\ldots k-1, (i_{k}=i_1), \forall n\geq 0\right]>0.
\end{equation}
Since $(\tau^{\varepsilon}_{n})_{n\geq 1}$ is monotone and bounded, then we can define $\tau:= \lim_{n\rightarrow\infty}\tau^{\varepsilon}_{n}$. Taking the limit with respect to $n$ in (\ref{limite}), we obtain
\begin{equation}\label{22b}
\P\left[Y^{i_l}_{\tau-}\leq Y^{i_{l+1}}_{\tau-}-g_{i_l,i_{l+1}}(\tau-), l=1,\ldots k-1, (i_{k}=i_1)\right]>0.
\end{equation}
But it is easy to verify that 
$$\left\{Y^{i_l}_{\tau-}\leq Y^{i_{l+1}}_{\tau-}-g_{i_l,i_{l+1}}(\tau-), l=1,\ldots k-1, (i_{k}=i_1)\right\}\subseteq\Big\{g_{i_1,i_2}(\tau-)+\cdots+g_{i_{k-1},i_1}(\tau-)\leq 0\Big\},$$
then from (\ref{22b}) we have
\begin{equation}
\P\left[g_{i_1,i_2}(\tau-)+\cdots+g_{i_k,i_1}(\tau-)\leq0\right]>0.\notag
\end{equation}
Since $g_{ij}\geq\gamma>0$ $\mathbb{P}$-a.s., we have a contradiction. Therefore, $\mathcal{S}^{\varepsilon}$ is admissible.$\hspace{4cm}\blacksquare$

Our next result has to do with a so-called verification theorem for the switching problem (\ref{pcg}) in the context of c\`adl\`ag cost functions
\begin{THM}\label{1.2}
The $q$ $\mathscr{S}^{2}$-processes $(Y^i_{\cdot}:=\left(Y^{i}_{t}\right)_{t\leq T}, i=1,\ldots,q)$ in Theorem \ref{te} are unique and they have the following relation with the switching problem (\ref{pcg}):
\begin{enumerate}[label=(\roman*)]
\item For each $i\in \mathcal{I}$,  
\begin{equation}\label{21}
Y^{i}_{0}=\sup_{\mathcal{S}\in\mathcal{A}_i}J(\mathcal{S},i). 
\end{equation}
\item The $\varepsilon$-strategy $\mathcal{S}^{\varepsilon}$ defined in Proposition \ref{admissible} forms an \textbf{$\varepsilon$}-optimal strategy, i.e., for $\mathcal{S}^{\varepsilon}=(\tau^{\varepsilon}_{n},\xi^{\varepsilon}_{n})_{n\geq 0}$, 
\begin{equation}
J\left(\mathcal{S}^{\varepsilon},i\right)\geq \sup\limits_{\mathcal{S}\in\mathcal{A}_i}J(\mathcal{S},i)-\varepsilon.
\end{equation}
\end{enumerate}
\end{THM}
\emph{Proof.}
\begin{enumerate}[label=(\roman*)]
\item Assuming that at time $t=0$ the system is in mode $i$, it follows by (\ref{a}) that, for any $0\leq t\leq T$,
\begin{equation}
Y^{i}_{t}+\int^{t}_{0}\psi_i(s)ds=\esssup_{\tau\geq t}\expect\bigg[\displaystyle\int^{\tau}_{0}\psi_i\left(s\right)ds+\max_{k\in\mathcal{I}^{-i}}\left(Y^{k}_{\tau}-g_{ik}\left(\tau\right)\right)\mathbf{1}_{\left[\tau<T\right]}\Big|\mathcal{F}_t\bigg].\notag
\end{equation}
Since $Y^{i}_{0}$ is $\mathcal{F}_0$-measurable, it is a $\mathbb{P}$-a.s. constant, that is, $Y^{i}_{0}=\mathbb{E}\left[Y^{i}_{0}\right]$. Now take $\mathcal{S}^{\varepsilon}$ defined in Proposition \ref{admissible}. Arguing similarly to Theorem \ref{te}, we can deduce 
\begin{equation}\label{b5}
\begin{array}{ll}
Y^{i}_{0} & \leq \expect\bigg[\displaystyle\int^{\tau^{\varepsilon}_{1}}_{0}\psi_i\left(s\right)ds+\max_{k\in\mathcal{I}^{-i}}\left(Y^{k}_{\tau^{\varepsilon}_{1}}-g_{ik}\left(\tau^{\varepsilon}_{1}\right)\right)\mathbf{1}_{\left[\tau^{\varepsilon}_{1}<T\right]}\bigg]+\frac{\varepsilon}{2} \\ \\ & = \expect\bigg[\displaystyle\int^{\tau^{\varepsilon}_{1}}_{0}\psi_i\left(s\right)ds+\left(Y^{\xi^{\varepsilon}_{1}}_{\tau^{\varepsilon}_{1}}-g_{i\xi^{\varepsilon}_{1}}\left(\tau^{\varepsilon}_{1}\right)\right)\mathbf{1}_{\left[\tau^{\varepsilon}_{1}<T\right]}\bigg]+\frac{\varepsilon}{2}.
\end{array}
\end{equation}
The rest of the proof uses the same arguments as in the proof of Theorem \ref{te}. Namely, for every $\tau^{\varepsilon}_{1} \leq t \leq T$, we can deduce
\begin{equation}
Y^{\xi^{\varepsilon}_{1}}_{t}=\esssup_{\tau\geq t}\expect\bigg[\displaystyle\int^{\tau}_{t}\psi_{\xi^{\varepsilon}_{1}}\left(s\right)ds+\max_{j\in\mathcal{I}^{-\xi^{\varepsilon}_{1}}}\left(Y^{j}_{\tau}-g_{\xi^{\varepsilon}_{1}j}\left(\tau\right)\right)\mathbf{1}_{\left[\tau<T\right]}\Big|\mathcal{F}_t\bigg].\notag
\end{equation}
Then, from the definition of $\tau^{\varepsilon}_{2}$ and since  $(Y^{\xi^{\varepsilon}_{1}}_{t} + \int^{t}_{\tau^{\varepsilon}_{1}}\psi_{\xi^{\varepsilon}_{1}}\left(s\right)ds)_{\tau^{\varepsilon}_{1}\leq t \leq\tau^{\varepsilon}_{2}}$ is a martingale, we get
\begin{equation}\label{b6}
\begin{array}{ll}
Y^{\xi^{\varepsilon}_{1}}_{\tau^{\varepsilon}_{1}} & \leq \expect\bigg[\displaystyle\int^{\tau^{\varepsilon}_{2}}_{\tau^{\varepsilon}_{1}}\psi_{\xi^{\varepsilon}_{1}}\left(s\right)ds+\max_{j\in\mathcal{I}^{-\xi^{\varepsilon}_{1}}}\left(Y^{j}_{\tau^{\varepsilon}_{2}}-g_{\xi^{\varepsilon}_{1}j}\left(\tau^{\varepsilon}_{2}\right)\right)\mathbf{1}_{\left[\tau^{\varepsilon}_{2}<T\right]}\Big|\mathcal{F}_{\tau^{\varepsilon}_{1}}\bigg]+\frac{\varepsilon}{4} \\ \\ & =\expect\bigg[\displaystyle\int^{\tau^{\varepsilon}_{2}}_{\tau^{\varepsilon}_{1}}\psi_{\xi^{\varepsilon}_{1}}\left(s\right)ds+\left(Y^{\xi^{\varepsilon}_{2}}_{\tau^{\varepsilon}_{2}}-g_{\xi^{\varepsilon}_{1}\xi^{\varepsilon}_{2}}\left(\tau^{\varepsilon}_{2}\right)\right)\mathbf{1}_{\left[\tau^{\varepsilon}_{2}<T\right]}\Big|\mathcal{F}_{\tau^{\varepsilon}_{1}}\bigg]+\frac{\varepsilon}{4}.
\end{array}
\end{equation}
Plugging (\ref{b6}) into (\ref{b5}) and noting that $\mathbf{1}_{\left[\tau^{\varepsilon}_{1}<T\right]}$ is $\mathcal{F}_{\tau^{\varepsilon}_{1}}$-measurable, it follows that: 
\begin{equation}
\begin{array}{ll}
Y^{i}_{0} & \leq \expect\bigg[\displaystyle\int^{\tau^{\varepsilon}_{1}}_{0}\psi_i\left(s\right)ds-g_{i\xi^{\varepsilon}_{1}}\left(\tau^{\varepsilon}_{1}\right)\mathbf{1}_{\left[\tau^{\varepsilon}_{1}<T\right]}\bigg] \\ \\ & \hspace{1cm} +\expect\bigg[\displaystyle\int^{\tau^{\varepsilon}_{2}}_{\tau^{\varepsilon}_{1}}\psi_{\xi^{\varepsilon}_{1}}\left(s\right)ds+\left(Y^{\xi^{\varepsilon}_{2}}_{\tau^{\varepsilon}_{2}}-g_{\xi^{\varepsilon}_{1}\xi^{\varepsilon}_{2}}\left(\tau^{\varepsilon}_{2}\right)\right)\mathbf{1}_{\left[\tau^{\varepsilon}_{2}<T\right]}\bigg]+\varepsilon\left(\frac{1}{2}+\frac{1}{4}\right). \\ \\ & = \expect\Bigg[\displaystyle\sum^{1}_{j=0}\bigg(\int^{\tau^{\varepsilon}_{j+1}}_{\tau^{\varepsilon}_{j}}\psi_{\xi^{\varepsilon}_{j}}(s)ds-g_{\xi^{\varepsilon}_{j}\xi^{\varepsilon}_{j+1}}\left(\tau^{\varepsilon}_{j}\right)\mathbf{1}_{\left[\tau^{\varepsilon}_{j+1}<T\right]}\bigg)+Y^{\xi^{\varepsilon}_{2}}_{\tau^{\varepsilon}_{2}}\mathbf{1}_{\left[\tau^{\varepsilon}_{2}<T\right]}\Bigg]+\varepsilon\left(\frac{1}{2}+\frac{1}{4}\right)
\end{array}\notag
\end{equation}
since $[\tau^{\varepsilon}_{2}<T]\subset [\tau^{\varepsilon}_{1}<T]$. Repeating this procedure $n$ times, we obtain
\begin{equation}
Y^{i}_{0}\leq \expect\Bigg[\displaystyle\sum^{n-1}_{j=0}\bigg(\int^{\tau^{\varepsilon}_{j+1}}_{\tau^{\varepsilon}_{j}}\psi_{\xi^{\varepsilon}_{j}}\left(s\right)ds-g_{\xi^{\varepsilon}_{j}\xi^{\varepsilon}_{j+1}}(\tau^{\varepsilon}_{j+1})\mathbf{1}_{\left[\tau^{\varepsilon}_{j+1}<T\right]}\bigg)+Y^{\xi^{\varepsilon}_{n}}_{\tau^{\varepsilon}_{n}}\mathbf{1}_{\left[{\tau^{\varepsilon}_{n}<T}\right]}\Bigg]+\varepsilon\left(\frac{1}{2}+\cdots +\frac{1}{2^{n}}\right).\notag
\end{equation}
Taking liminf as $n\rightarrow\infty$ we obtain
\begin{equation}\label{25a}
Y^{i}_{0}\leq \expect\Bigg[\displaystyle\sum^{\infty}_{j=0}\bigg(\int^{\tau^{\varepsilon}_{j+1}\wedge T}_{\tau^{\varepsilon}_{j}\wedge T}\psi_{\xi^{\varepsilon}_{j}}\left(s\right)ds-g_{\xi^{\varepsilon}_{j}\xi^{\varepsilon}_{j+1}}(\tau^{\varepsilon}_{j+1})\mathbf{1}_{\left[\tau^{\varepsilon}_{j+1}<T\right]}\bigg)\bigg]+\varepsilon.
\end{equation}
By Proposition \ref{admissible} we can take supremum over all admissible strategies $\mathcal{A}_i$, to obtain  
\begin{equation}
\begin{array}{ll}
Y^{i}_{0} & \leq \sup\limits_{\mathcal{S}\in \mathcal{A}_i}\expect\Bigg[\displaystyle\sum^{\infty}_{j=0}\bigg(\int^{\tau_{j+1}\wedge T}_{\tau_{j}\wedge T}\psi_{\xi_{j}}\left(s\right)ds-g_{{\xi_j}\xi_{j+1}}(\tau_{j+1})\mathbf{1}_{\left[\tau_{j+1}<T\right]}\bigg)\Bigg]+\varepsilon \\  \\ & =\sup\limits_{\mathcal{S}\in\mathcal{A}_i}J(\mathcal{S},i)+\varepsilon.\notag
\end{array}
\end{equation}
Letting $\varepsilon\rightarrow 0$, it follows that $Y^{i}_{0}\leq\sup_{\mathcal{S}\in \mathcal{A}_i}J(\mathcal{S},i)$. The inverse inequality is analogous to the previous Theorem \ref{te}. Hence, the result follows.
\item From part (i), specifically, (\ref{21}) and inequality (\ref{25a}), we deduce 
\begin{equation}
\sup_{\mathcal{\mathcal{S}\in A}_i}J(\mathcal{S},i)-\varepsilon\leq J(\mathcal{S}^{\varepsilon},i)\leq \sup_{\mathcal{S}\in \mathcal{A}_i}J(\mathcal{S},i),\notag
\end{equation}
which proves (ii). $\hspace{14cm}\blacksquare$
\end{enumerate}
\section{Reflected Backward Stochastic Differential Systems}\label{sec3}
In this section we will provide the existence as well as 
uniqueness of the solution of the system of reflected backward stochastic differential equations (RBSDEs) of type 
\begin{equation}\label{RBSDEM5.4}
  \hspace{-0.2cm}	\begin{cases}
	\forall i\in\mathcal{I},\hspace{.2cm} \text{find}\hspace{.2cm} \left(Y^{i}_{\cdot},Z^{i}_{\cdot},K^{i}_{\cdot}\right)\hspace{.2cm} \text{such that}: &\\ Y^{i}_{\cdot},K^{i}_{\cdot}\in\mathscr{S}^{2}\quad \text{and}\quad Z^{i}_{\cdot}\in \mathcal{H}^{2,d};\hspace{0.1cm} K^{i}_{\cdot}\hspace{0.2cm} \text{is non-decreasing and}\quad K^{i}_{0}=0,& \\
    Y^{i}_{s}=h_i(X_T)+\displaystyle\int^{T}_{s}f_i(r,X_r,Y^{1}_{r},\ldots,Y^{q}_{r},Z^{i}_{r})dr+K^{i}_{T}-K^{i}_{s}-\int^{T}_{s}Z^{i}_{r}dB_r\hspace{0.4cm} \text{for all }  0\leq s\leq T,& \\ 
		Y^{i}_{s}\geq\max\limits_{k\in\mathcal{I}^{-i}}\left\{Y^{k}_{s}-\gamma_{ik}(s,X_{s})\right\}\hspace{1cm} \text{for all }  0\leq s\leq T, & \\ \text{and if}\hspace{0.2cm} K^{i}_{\cdot}=K^{i,c}_{\cdot}+K^{i,d}_{\cdot},\hspace{0.1cm}\text{where}\hspace{0.2cm}K^{i,c}_{\cdot}\hspace{0.1cm}(\text{resp.}\hspace{0.2cm}K^{i,d}_{\cdot})\hspace{0.2cm}\text{is the continuous} & \\ \text{(resp. purely discontinuous) part of}\hspace{0.2cm} K^{i}_{\cdot},\hspace{0.2cm} \text{then:} \\
		\displaystyle\int^{T}_{0}\Big(Y^{i}_{r}-\max_{k\in\mathcal{I}^{-i}}\big\{Y^{k}_{r}-\gamma_{ik}(r,X_{r})\big\}\Big)dK^{i,c}_{r}=0. & \\
		 \Delta_sY_{\cdot}:=Y_s-Y_{s^{-}}=-\Big(\max\limits_{k\in\mathcal{I}^{-i}}\left\{Y^{k}_{s}-\gamma_{ik}(s,X_{s})\right\}-Y^{i}_s\Big)^{+} \hspace{1cm} \text{for all }  0\leq s\leq T,
  \end{cases}
\end{equation}
in which the associated barriers are \emph{c\`adl\`ag} processes. This system is connected with the previous switching problem. Actually when $(f_i)_{\ici}$ do not depend on $(Y^i)_{\ici}$, the system \eqref{RBSDEM5.4} is exactly the translation of the verification Theorem \ref{te} in terms of reflected BSDEs as it is well-known that the Snell envelope can be expressed through reflected BSDEs (see e.g. El Karoui \cite{karoui1981aspects} or Hamad\`ene \cite{hamadene2002reflected}). On the other hand, this form of system \eqref{RBSDEM5.4} allows to consider switching problems when the cost functions are of risk sensitive type  (utility functions) ---see El Karoui and Hamad\`ene \cite{el2003bsdes}.
\medskip

\noindent To begin with our analysis, we will first introduce the following assumptions relate to the items involved in (\ref{RBSDEM5.4}):
\medskip

\begin{assumptionH}\label{assH} \
\begin{enumerate}[label=(H{{\arabic*}})]
\item\label{H1}: The stochastic process $(X_t)_{t\geq 0}$ is in $\mathcal{H}^{2,r}$ for any $r\in \mathbb{N}$.  
\item\label{H2}: For any $i\in\mathcal{I}$, the function $f_i:[0,T]\times\mathbb{R}^{r}\times\mathbb{R}^{q}\times\mathbb{R}^{d}\rightarrow\mathbb{R}$ satisfies:
	\begin{enumerate}[label=(\roman*)]
\item\label{H2i}$\left(t,x\right)\mapsto f_i\left(t,x,y^{1},\ldots,y^{q},z\right)$ is continuous uniformly with respect to $(y^{1},\ldots,y^{q},z)$; 
\item\label{H2ii} $f_i$ is uniformly Lipschitz continuous with respect to $(y^{1},\ldots,y^{q},z)$, i.e., for some $C\geq 0$,
\begin{equation}
\left|f_i(t,x,y^{1},\ldots,y^{q},z)-f_i(t,x,\bar{y}^{1},\ldots,\bar{y}^{q},\bar{z})\right| \leq  C\left(\left|y^{1}-\bar{y}^{1}\right|+\cdots+\left|y^{q}-\bar{y}^{q}\right|+\left|z-\bar{z}\right|\right). \notag
\end{equation}
\item\label{H2iii} the mapping $(t,x)\mapsto f_i(t,x,0,\ldots,0)$ is Borel measurable and of polynomial growth.
\item\label{H2iv} \emph{Monotonicity}: For all $\ici$, for all $k\in\mathcal{I}^{-i}$, the mapping $y_k\mapsto f_i(t,x, y_1,\ldots, y_{k-1},y_k,y_{k+1},\ldots, y_q,z)$ is non-decreasing whenever the other components $(t, x, y_1,\ldots,y_{k-1},y_{k+1},\ldots, y_q,z)$ are fixed.
\end{enumerate}
\item\label{H3}: For each $i,k\in\mathcal{I}$, the function $\gamma_{ik}:[0,T]\times\mathbb{R}^{r}\rightarrow\mathbb{R}$ is bounded from below, i.e. there exists a real constant $\gamma>0$ such that, $\gamma_{ik}\geq \gamma$. Furthermore it is c\`adl\`ag in $t$, continuous and of polynomial growth in $x$.
\item\label{H4}: For each $\ici$, the function $h_i:\mathbb{R}^{r}\rightarrow\mathbb{R}$ is continuous with polynomial growth and satisfies
	\begin{equation}
	\forall x\in\mathbb{R}^{r},h_i(x)\geq \max_{k\in\mathcal{I}^{-i}}\left(h_k(x)-\gamma_{ik}(T,x)\right). \notag \qed
	\end{equation}
\end{enumerate}
\end{assumptionH}

Note that in the \eqref{RBSDEM5.4} the process $X$ does not play a specific role. We consider this form of system \eqref{RBSDEM5.4} only in the perspective to deal with the Hamilton-Jacobi-Bellman system associated with the switching problem.
\begin{PROP}\label{5.1}
Under Assumptions (\ref{assH}), the system of RBSDEs (\ref{RBSDEM5.4}) has a solution $(Y^{i}_{\cdot},Z^{i}_{\cdot},K^{i}_{\cdot} )_{i=1,\ldots,q}$.
\end{PROP}
\noindent \emph{Proof.} To begin with, we first consider the following standard BSDEs:
	\begin{equation}\label{BSDEMS}
  	\hspace{-0.1cm}\begin{cases}
	  (\overline{Y}_{\cdot},\overline{Z}_{\cdot})\in \mathscr{S}^{2}\times \mathcal{H}^{2,d}\sp ;& \\
    \overline{Y}_{s}=\max\limits_{i=1,\ldots,q}h_i(X_T)+\displaystyle\int^{T}_{s}\big[\max_{i=1,\ldots,q}f_i\big](r,X_{r},\overline{Y}_{r},\ldots,\overline{Y}_{r},\overline{Z}_{r})dr-\int^{T}_{s}\overline{Z}_{r}dB_r, \hspace{0.25cm}\text{for all } s\leq T,
  \end{cases}
\end{equation}
and
	\begin{equation}\label{BSDEMM}
  	\hspace{-0.2cm}\begin{cases}
	(\underline{Y}_{\cdot},\underline{Z}_{\cdot})\in \mathscr{S}^{2}\times \mathcal{H}^{2,d}\sp ;
& \\
    \underline{Y}_{s}=\min\limits_{i=1,\ldots,q}h_i(X_T)+\displaystyle\int^{T}_{s}\big[\min_{i=1,\ldots,q}f_i\big](r,X_{r},\underline{Y}_{r},\ldots,\underline{Y}_{r},\underline{Z}_{r})dr-\int^{T}_{s}\underline{Z}_{r}dB_r, \hspace{0.25cm}\text{for all } s\leq T.
  \end{cases}
\end{equation}
It is easy to verify that under (\ref{assH}) the data of \eqref{BSDEMS} and  \eqref{BSDEMS} satisfy the conditions of Pardoux and Peng's result \cite{pardoux1990adapted}, pp. 59-60 and then in virtue of Theorem 4.1 of this same reference, we claim the existence and uniqueness of solutions of both (\ref{BSDEMS}) and (\ref{BSDEMM}). To solve the system (\ref{RBSDEM5.4}), we shall use an iterative method and regard (\ref{RBSDEM5.4}) as a limit system. To this end, for any $i\in\mathcal{I}$, we set $Y^{i,0}_{\cdot}:=\underline{Y}_{\cdot}$, and for $n\geq 1$, we seek a triplet $(Y^{i,n}_{\cdot},Z^{i,n}_{\cdot},K^{i,n}_{\cdot})$ such that, for $n\geq 1$
	\begin{equation}
  	\hspace{-0.1cm}\begin{cases}
	  Y^{i,n}_{\cdot},K^{i,n}_{\cdot}\in\mathscr{S}^{2}\quad \text{and}\quad Z^{i,n}_{\cdot}\in \mathcal{H}^{2,d};\hspace{0.1cm} K^{i,n}_{\cdot}\hspace{0.2cm} \text{is non-decreasing with}\quad K^{i}_{0}=0,& \\
    Y^{i,n}_{s}=h_i(X_T)+\displaystyle\int^{T}_{s}f_i(r,X_{r},Y^{1,n-1}_{r},\ldots,Y^{i-1,n-1}_{r},Y^{i,n}_{r},Y^{i+1,n-1}_{r},\ldots Y^{q,n-1}_{r},Z^{i,n}_{r})dr& \\ \hspace{1.1cm}+K^{i,n}_{T}-K^{i,n}_{s}-\int^{T}_{s}Z^{i,n}_{r}dB_r,\hspace{0.2cm}\text{for all } 0\leq s\leq T; & \\ 
		Y^{i,n}_{s}\geq\max\limits_{k\in\mathcal{I}^{-i}}\left\{Y^{k,n-1}_{s}-\gamma_{ik}(s,X_{s})\right\} , \hspace{0.25cm}\text{for all } 0\leq s\leq T, & \\ \text{and if}\hspace{0.2cm} K^{i,n}_{\cdot}=K^{i,n,c}_{\cdot}+K^{i,n,d}_{\cdot},\hspace{0.1cm}\text{where}\hspace{0.2cm}K^{i,n,c}_{\cdot}\hspace{0.1cm}(\text{resp.}\hspace{0.2cm} K^{i,n,d}_{\cdot})\hspace{0.2cm}\text{is the continuous} & \\ \text{(resp. purely discontinuous) part of}\hspace{0.2cm} K^{i,n}_{\cdot},\hspace{0.2cm} \text{then:}& \\
		\displaystyle\int^{T}_{0}\left(Y^{i,n}_{r}-\max\limits_{k\in\mathcal{I}^{-i}}\left\{Y^{k,n-1}_{r}-\gamma_{ik}(r,X_{r})\right\}\right)dK^{i,n,c}_{r}=0 \,; & \\
		 \Delta_sY_{\cdot}:=Y^{i,n}_s-Y^{i,n}_{s^{-}}=-\Big(\max\limits_{k\in\mathcal{I}^{-i}}\left\{Y^{k,n-1}_{s}-\gamma_{ik}(s,X_{s})\right\}-Y^{i,n}_s\Big)^{+}, \hspace{0.25cm}\text{for all} 0\leq s\leq T. \notag
  \end{cases}
\end{equation}
Note that for each $k\in\mathcal{I}$ the process $Y^{k,0}_{\cdot}$ is given. Then, by letting 
$$\tilde{f}_i(s,Y^{i,1}_s,Z^{i,1}_s):=f_i(s,X_{s},Y^{1,0}_{s},\ldots,Y^{i-1,0}_{s},Y^{i,1}_{s},Y^{i+1,0}_{s},\ldots Y^{q,0}_{s},Z^{i,1}_{s})$$ for $i\in\mathcal{I}$,
the data of the RBSDE associated with $(Y^{i,1}_{\cdot},Z^{i,1}_{\cdot},K^{i,1}_{\cdot})$ satisfy the assumptions in Hamad\`ene \cite{hamadene2002reflected}, Theorem 1.4, and hence the processes $(Y^{i,1}_{\cdot},Z^{i,1}_{\cdot},K^{i,1}_{\cdot})$ do exist. Next, using the comparison theorem of solutions of BSDEs (see e.g. Theorem 2.2 in El Karoui et al. \cite{karoui1997backward}) we deduce that for any $i\in \mathcal{I}$, $Y^{i,0}_{\cdot}\leq Y^{i,1}_{\cdot}$. Besides, as $f_i$ satisfies the monotonicity property \ref{H2}-\ref{H2iv} and using again the comparison of solutions of RBSDEs (see Theorem 1.5 in Hamad\`ene \cite{hamadene2002reflected}) we obtain by induction that:
\begin{equation}\label{increasing}
\text{for all}\hspace{0.2cm} n\geq 0 \hspace{0.2cm}\text{and}\hspace{0.2cm}i\in\mathcal{I}, \hspace{0.4cm} Y^{i,n}_{\cdot}\leq Y^{i,n+1}_{\cdot}.
\end{equation}
On the other hand, the process $(\overline{Y}_{\cdot},\overline{Z}_{\cdot})$ in (\ref{BSDEMS}), can be regarded as the triplet $((\overline{Y}^{i}_{\cdot},\overline{Z}^{i}_{\cdot},0))_{i\in\mathcal{I}}$ (i.e. $K^{i}_{\cdot}=0$), which is solution for the system of RBSDEs with data $$([\max\limits_{i=1,\ldots,q}f_i](t,X_{t},y^{1},\ldots,y^{q},z),\gamma_{ik}(t,X_{t}),\max\limits_{i=1,\ldots,q}h_i(X_T)),\quad i,k\in \mathcal{I}.$$ Note that
\begin{equation}
\begin{array}{ll}
f_i(t,\overline{Y}_t,\overline{Z}_t)& :=f_i(t,X_t,Y^{1,0}_{t},\ldots,\overline{Y}^{i}_t,\ldots,Y^{q,0}_{t},\overline{Z}_{t})  \\ & \leq \max_{i=1,\ldots,q}f_i(t,X_t,\overline{Y}^{1}_t,\ldots,\overline{Y}^{i}_{t},\ldots,\overline{Y}^{q}_{t},\overline{Z}_{t}) \\ &
:= \max_{i=1,\ldots,q}f_i(t,\overline{Y}_t,\overline{Z}_t),
\end{array}\notag
\end{equation}
since $f_i$ satisfies the monotonicity property \ref{H2}-\ref{H2iv} and due that, for each $k\in\mathcal{I}^{-i}$ (the fixed processes), $Y^{k,0}_{\cdot}=\underline{Y}_{\cdot}\leq \overline{Y}_{\cdot}$. Therefore, by comparison  Theorem 1.5 in Hamad\`ene \cite{hamadene2002reflected}, we get that $Y^{i,1}_{\cdot}\leq \overline{Y}_{\cdot}$. In general, through an induction procedure, we can obtain for all $n\geq 0$ and $i\in\mathcal{I}$, $Y^{i,n}_{\cdot}\leq \overline{Y}_{\cdot}$ and hence
\begin{equation}\label{growth}
\underline{Y}_{\cdot}=Y^{i,0}_{\cdot}\leq Y^{i,n}_{\cdot}\leq Y^{i,n+1}_{\cdot}\leq \overline{Y}_{\cdot}.
\end{equation}
Arguing as in Theorem \ref{te}, we can see that there exists $Y^{i}_{\cdot}$ such that $Y^{i,n}_{\cdot}\nearrow Y^{i}_{\cdot}$ and $\expect\big[\sup_{0\leq t\leq T}\left|Y^{i}_{t}\right|^{2}\big]<\infty$. Therefore, 
using Peng's monotonic limit theorem (see Theorem 2.1 and Theorem 3.6 in Peng \cite{peng1999monotonic}), we deduce that for any $i\in\mathcal{I}$, the limit process $Y^{i}_{\cdot}$  is a c\`adl\`ag process and there exists $(Z^{i}_{\cdot},K^{i}_{\cdot})\in\mathcal{H}^{2,d}\times\mathscr{S}^{2}$ with $K^{i}$ non-decreasing process and $K^{i}_{0}=0$ such that: $\forall s\le T$, 
\begin{equation}
  	\hspace{-0.03cm}\begin{cases}
Y^{i}_{s}=h_i(X_T)+\displaystyle\int^{T}_{s}f_i(r,X_{r},Y^{1}_{r},\ldots, Y^{i}_{r},\ldots,Y^{q}_{r},Z^{i}_{r})dr+K^{i}_{T}-K^{i}_{s}-\int^{T}_{s}Z^{i}_{r}dB_r\,;& \\ 
		Y^{i}_{s}\geq\max\limits_{k\in\mathcal{I}^{-i}}\left\{Y^{k}_{s}-\gamma_{ik}(s,X_{s})\right\}.
  \end{cases}\notag
\end{equation}
Now we claim that $\left(Y^{i},Z^{i},K^{i}\right)_{i=1,\ldots,q}$ is, in fact, the desired solution of (\ref{RBSDEM5.4}). Indeed, consider the RBSDEs at the $i$-th variable and the other variables $Y^{1},\dots,Y^{i-1},Y^{i+1},\dots,Y^{q}$ fixed, that is to say 
	\begin{equation}\label{rbsdefija}
  \begin{cases}
	\forall i\in\mathcal{I}, \hspace{.2cm} \text{find}\hspace{.2cm} \left(\tilde{Y}^{i}_{\cdot},\tilde{Z}^{i}_{\cdot},\tilde{K}^{i}_{\cdot}\right)\hspace{.2cm}\text{such that}: \\ \tilde{Y}^{i}_{\cdot},\tilde{K}^{i}_{\cdot}\in\mathscr{S}^{2}\quad \text{and}\quad \tilde{Z}^{i}_{\cdot}\in \mathcal{H}^{2,d}; \tilde{K}^{i}_{\cdot}\hspace{.2cm} \text{is non-decreasing and} \hspace{.2cm}\tilde{K}^{i}_{0}=0 & \\
    \tilde{Y}^{i}_{s}=h_i(X_T)+\displaystyle\int^{T}_{s}f_i(r,X_{r},Y^{1}_{r},\ldots,Y^{i-1}_{r},\tilde{Y}^{i}_{r},Y^{i+1}_{r},\ldots Y^{q}_{r},Z^{i}_{r})dr+\tilde{K}^{i}_{T}-\tilde{K}^{i}_{s}-\int^{T}_{s}\tilde{Z}^{i}_{r}dB_r,\,\,\text{for all } 0\leq s\leq T;& \\ 
		\tilde{Y}^{i}_{s}\geq\max\limits_{k\in\mathcal{I}^{-i}}\left\{Y^{k}_{s}-\gamma_{ik}(s,X_{s})\right\}\hspace{0.5cm}\text{for all } 0\leq s\leq T  \\ \text{and if}\hspace{0.2cm} \tilde{K}^{i}_{\cdot}=\tilde{K}^{i,c}_{\cdot}+K^{i,d}_{\cdot},\hspace{0.1cm}\text{where}\hspace{0.2cm}\tilde{K}^{i,c}_{\cdot}\hspace{0.1cm}(\text{resp.} \tilde{K}^{i,d}_{\cdot})\hspace{0.2cm}\text{is the continuous} & \\ \text{(resp. purely discontinuous) part of}\hspace{0.2cm} \tilde{K}^{i},\hspace{0.2cm} \text{then:} \\
		\displaystyle\int^{T}_{0}\Big(\tilde{Y}^{i}_{r}-\max\limits_{k\in\mathcal{I}^{-i}}\big\{Y^{k}_{r}-\gamma_{ik}(r,X_{r})\big\}\Big)d\tilde{K}^{i,c}_{r}=0. \\
		 \Delta_s\tilde{Y}_{\cdot}:=\tilde{Y}^{i}_s-\tilde{Y}^{i}_{s^{-}}=-\Big(\max\limits_{k\in\mathcal{I}^{-i}}\left\{Y^{k}_{s}-\gamma_{ik}(s,X_{s})\right\}-\tilde{Y}^{i}_s\Big)^{+}\hspace{0.5cm}\text{for all } 0\leq s\leq T. 
		 \end{cases}
\end{equation} 
The solution of $(\ref{rbsdefija})$ do exist by using again Theorem 1.4 in Hamad\`ene \cite{hamadene2002reflected}. Such a solution $\tilde{Y}^{i}_{\cdot}$ becomes the smallest $f_i$-supermartingale that dominates $\max_{k\in\mathcal{I}^{-i}}\left\{Y^{k}_{s}-\gamma_{ik}(s,X_{s})\right\}$ (for more details on this last assertion, see Peng and Xu \cite{peng2005smallest}). Whence $\tilde{Y}^{i}_{t}\leq Y^{i}_{t}$. On the other hand, since $Y^{i,n-1}_{t}\leq Y^{i}_{t}$ for any $i\in \mathcal{I}$ and $n\geq 1$, we get
    \begin{equation}
  \max\limits_{k\in\mathcal{I}^{-i}}\big\{Y^{k,n-1}_{s}-\gamma_{ik}(s,X_{s})\big\}\leq \max\limits_{k\in\mathcal{I}^{-i}}\big\{Y^{k}_{s}-\gamma_{ik}(s,X_{s})\big\}.\notag
\end{equation}
Also observe that assumptions (\ref{H2})-(\ref{H2iv}) yields that
$$f_i(t,x,Y^{1,n-1}_{t},\ldots,\tilde{Y}^{i}_{t},\ldots,Y^{q,n-1}_{t},Z^{i}_{t})\leq f_i(t,x,Y^{1}_{t},\ldots,\tilde{Y}^{i}_{t},\ldots,Y^{q}_{t},Z^{i}_{t}).$$
Then using again the comparison theorem for RBSDEs given in Theorem 1.5 in Hamad\`ene \cite{hamadene2002reflected}, we have $Y^{i,n}_{t}\leq \tilde{Y}^{i}_{t}$. This implies that $Y^{i}_{t}\leq \tilde{Y}^{i}_{t}$, and hence $\tilde{Y}^{i}_{t}=Y^{i}_{t}$. Moreover, this also implies that $\tilde{Z}^{i}_{t}=Z^{i}_{t}$ and $\tilde{K}^{i}_{t}=K^{i}_{t}$ for any $0\leq t\leq T$, $\mathbb{P}$-a.s. This proves the existence of solution for (\ref{RBSDEM5.4}). $\qed$
\medskip

We now provide a representation result for the solutions of system (\ref{RBSDEM5.4}) and, as a by-product, we obtain the uniqueness. For later use, let us fix $\textbf{u}_{.}:=\left(u^{1}_{.},\ldots,u^{q}_{.}\right)$ in $\mathcal{H}^{2,q}$ and let us consider the following system of RBSDEs: 
\begin{equation}\label{5.9}
  	\begin{cases}
		\forall i\in\mathcal{I}, \hspace{.2cm} \text{find}\hspace{.2cm} \left(Y^{\textbf{u},i}_{\cdot},Z^{\textbf{u},i}_{\cdot},K^{\textbf{u},i}_{\cdot}\right)\in \mathscr{S}^2\times\mathscr{S}^2\times\mathcal{H}^{2,d}\hspace{.2cm}\text{such that}: & \\ Y^{\textbf{u},i}_{s}=h_i(X_T)+\displaystyle\int^{T}_{s}f_i(r,X_{r},\textbf{u}_r,Z^{\textbf{u},i}_{r})dr+K^{\textbf{u},i}_{T}-K^{\textbf{u},i}_{s}-\int^{T}_{s}Z^{\textbf{u},i}_{r}dB_r\hspace{0.5cm}\text{for all } 0\leq s\leq T;& \\ 
		Y^{\textbf{u},i}_{s}\geq\max\limits_{k\in\mathcal{I}^{-i}}\left\{Y^{\textbf{u},k}_{s}-\gamma_{ik}(s,X_{s})\right\}\hspace{0.5cm}\text{for all } 0\leq s\leq T. & \\ 
		\text{and if}\hspace{0.2cm} K^{\textbf{u},i}_{\cdot}=K^{\textbf{u},i,c}_{\cdot}+K^{\textbf{u},i,d}_{\cdot},\hspace{0.1cm}\text{where}\hspace{0.2cm}K^{i,\textbf{u},c}_{\cdot}\hspace{0.1cm}(\text{resp.} K^{i,\textbf{u},d}_{\cdot})\hspace{0.2cm}\text{is the continuous} & \\ \text{(resp. purely discontinuous) part of}\hspace{0.2cm} K^{\textbf{u},i}_{\cdot},\hspace{0.2cm} \text{then:}& \\
		\displaystyle\int^{T}_{0}\Big(Y^{\textbf{u},i}_{r}-\max\limits_{k\in\mathcal{I}^{-i}}\big\{Y^{\textbf{u},k}_{r}-\gamma_{ik}(r,X_{r})\big\}\Big)dK^{i,\textbf{u},c}_{r}=0. & \\
		 \Delta_s\tilde{Y}^{\textbf{u}}_{\cdot}:=\tilde{Y}^{\textbf{u},i}_s-\tilde{Y}^{\textbf{u},i}_{s^{-}}=-\Big(\max\limits_{k\in\mathcal{I}^{-i}}\left\{Y^{\textbf{u},k}_{s}-\gamma_{ik}(s,X_{s})\right\}-\tilde{Y}^{\textbf{u},i}_s\Big)^{+}\hspace{0.5cm}\text{for all } 0\leq s\leq T. 
  \end{cases}
\end{equation}
Observe that $f_i$ does not depend on $Y^{1},\ldots, Y^{q}$. Let $s\leq T$ be fixed, $i\in\mathcal{I}$ and let $\mathcal{D}^{i}_{s}$ be the following set of strategies as in (\ref{stra}), such that:
\begin{equation}
\mathcal{D}^{i}_{s}:=\left\{\alpha=\left(\theta_n,\kappa_n\right)_{n\geq 0}: \theta_0=s, \kappa_0=i \hspace{0.2cm}\text{and}\hspace{0.2cm} \mathbb{E}[(\mathtt{C}^{\alpha}_{T})^{2}]<\infty\right\}\notag
\end{equation}
where $\mathtt{C}^{\alpha}_{r}$, $r \leq T$, is the following cumulative costs up to time $r$, i.e.,
\begin{equation}
\mathtt{C}^{\alpha}_{r}:=\sum^{\infty}_{n=1}\gamma_{\kappa_{n-1},\kappa_n}(\theta_n,X_{\theta_n})\mathbf{1}_{[\theta_n\leq r]} \hspace{0.4cm} \text{for}\hspace{0.2cm} r<T\hspace{0.2cm} \text{and}\hspace{0.2cm} \mathtt{C}^{\alpha}_T=\lim_{r\rightarrow T} \mathtt{C}^{\alpha}_r, \hspace{0.2cm} \mathbb{P}\text{-a.s.}\notag
\end{equation}
Therefore and for any admissible strategy $\alpha\in\mathcal{D}^{i}_{s}$ we have:
\begin{equation}
\mathtt{C}^{\alpha}_{T}=\sum^{\infty}_{n=1}\gamma_{\kappa_{n-1},\kappa_n}(\theta_n,X_{\theta_n})\mathbf{1}_{[\theta_n<T]}.\notag
\end{equation}
Consider a strategy $\alpha=(\theta_n,\kappa_n)_{n\geq0}\in \mathcal{D}^{i}_{s}$ and let $(P^{\alpha}_{\cdot},Q^{\alpha}_{\cdot}):=(P^{\alpha}_{s} ,Q^{\alpha}_{s})_{s\leq T}$ be the solution of the following BSDE
\begin{equation}\label{repre1}
    \begin{cases}
        P^{\alpha}_{\cdot}\hspace{0.1cm} \text{is c\`adl\`ag and} \hspace{0.1cm}\expect\big[\sup_{s\leq T}\left|P^{\alpha}_{s}\right|^{2}\big]<\infty,\hspace{0.2cm} Q^{\alpha}_{\cdot}\in\mathcal{H}^{2,d};& \\
      P^{\alpha}_{s}=h_\alpha(X_T)+\displaystyle\int^{T}_{s}f_\alpha(r,X_{r},\textbf{u}_r,Q^{\alpha}_{r})dr-(\mathtt{C}^{\alpha}_{T}-\mathtt{C}^{\alpha}_{s})-\int^{T}_{s}Q^{\alpha}_{r}dB_r,\,s\le T,
  \end{cases}
\end{equation}
with
\begin{equation}\label{38a}
    \begin{array}{ll}
h_{\alpha}(x)= h_{\kappa_n}(x)\mathbf{1}_{[\theta_n< T\leq\theta_{n+1}]} \hspace{0.1cm}\text{and}\hspace{0.1cm} & \\ & \\ f_{\alpha}(r,x,v_1,\ldots,v_q,z):=\displaystyle\sum\limits^{\infty}_{n=0}f_{\kappa_n}(r,x,v_1,\ldots,v_q,z)\mathbf{1}_{[\theta_n\leq r<\theta_{n+1})}.
   \end{array}
\end{equation}
Making the change of variable $\bar{P}^{\alpha}_{\cdot}:= P^{\alpha}_{\cdot}-\mathtt{C}^{\alpha}_{\cdot}$, the equation in (\ref{repre1}) is transformed in a standard BSDE. Since $\mathtt{C}^{\alpha}$ is adapted and $\mathbb{E}[(\mathtt{C}^{\alpha}_{T})^{2}]<\infty$, we easily deduce the existence and uniqueness of the process $(P^{\alpha}_{\cdot},Q^{\alpha}_{\cdot})$. We then have the following representation for the solution of (\ref{5.9}).
\begin{PROP}\label{5.3}
Assume that for any $i, k\in\mathcal{I}$:
\begin{enumerate}[label=(\roman*)]
	\item $f_i$ satisfies \ref{H2}-\ref{H2ii},\ref{H2iii};
  \item $\gamma_{ik}$ (resp. $h_i$) satisfies \ref{H3} (resp. \ref{H4}).
\end{enumerate}
Then the solution of system of RBSDEs (\ref{5.9}) exists, it is unique and satisfies:
\begin{equation}\label{5.11}
 Y^{\textbf{u},i}_{s}=\esssup_{\alpha\in\mathcal{D}^{i}_{s}}\left\{P^{\alpha}_{s}-\mathtt{C}^{\alpha}_{s}\right\}\quad \forall s\leq T,\quad\forall i\in\mathcal{I}. 
\end{equation}
\end{PROP}
\emph{Proof.} Since $f_i$ does not depend on variables $Y^{1}_{\cdot},\ldots,Y^{q}_{\cdot}$, then, it trivially satisfies \ref{H2}-\ref{H2iv}. Then, by hypothesis $\emph{(i)}$ and $\emph{(ii)}$, and Proposition \ref{5.1}, the solution ($Y^{\textbf{u},i}_{\cdot},Z^{\textbf{u},i}_{\cdot},K^{\textbf{u},i}_{\cdot})$ of the system (\ref{5.9}) exists. Therefore, plugging an arbitrary strategy $\alpha\in\mathcal{D}^{i}_{s}$ in (\ref{5.9}), we obtain:
\begin{equation}\label{f1}
Y^{\textbf{u},i}_{s}\geq h_{\alpha}(X_{T})+\displaystyle\int^{T}_{s}f_{\alpha}(r,X_{r},\textbf{u}_{r},Z^{\alpha}_{r})dr+\tilde{K}^{\alpha}_{T}-\mathtt{C}^{\alpha}_{T}-\int^{T}_{s}Z^{\alpha}_{r}dB_r. 
\end{equation}
with $h_\alpha$ and $f_\alpha$ as in (\ref{38a}), and,
\begin{equation}
    \tilde{K}^{\alpha}_{T}=(K^{\textbf{u},i}_{\theta_1}-K^{\textbf{u},i}_{s})+\sum\limits^{\infty}_{n=1}(K^{\textbf{u},\kappa_n}_{\theta_{n+1}}-K^{\textbf{u},\kappa_n}_{\theta_n})\text{ and }Z^{\alpha}_{r}=\sum\limits^{\infty}_{n=0}Z^{\textbf{u},\kappa_n}_{r}\mathbf{1}_{[\theta_n\leq r<\theta_{n+1})}, \forall r\leq T.
\end{equation}
Adding $\mathtt{C}^{\alpha}_{s}$ from both sides of (\ref{f1}) and taking into account that $\tilde{K}^{\alpha}_{T}\geq 0$, we have 
\begin{equation}
\begin{array}{ll}
Y^{\textbf{u},i}_{s}+\mathtt{C}^{\alpha}_{s}\geq h_{\alpha}(X_{T})+\displaystyle\int^{T}_{s}f_{\alpha}(r,X_{r},\textbf{u}_{r},Z^{\alpha}_{r})dr-(\mathtt{C}^{\alpha}_{T}-\mathtt{C}^{\alpha}_{s})-\int^{T}_{s}Z^{\alpha}_{r}dB_r & \\ \hspace{1.7cm} = P^{\alpha}_{s}.
\end{array}\notag
\end{equation}
Therefore, we have
\begin{equation}\label{41a}
Y^{\textbf{u},i}_{s}\geq \esssup\limits_{\alpha\in\mathcal{D}^{i}_{s}}\left\{P^{\alpha}_{s}-\mathtt{C}^{\alpha}_{s}\right\},\hspace{0.3cm} \forall\alpha\in\mathcal{D}^{i}_{s}.
\end{equation}
Next let $\alpha^{\varepsilon}=(\theta^{\varepsilon}_{n},\kappa^{\varepsilon}_{n})_{n\geq 0}$ be the strategy defined recursively as follows (compare to the $\varepsilon$-strategy $\mathcal{S}^{\varepsilon}$ in Proposition \ref{admissible}): $\theta^{\varepsilon}_{0}=0,\kappa^{\varepsilon}_{0}= i$ and for $n\geq 0$,
\begin{equation}
\theta^{\varepsilon}_{n+1}=\inf\left\{s\geq\theta^{\varepsilon}_{n}: Y^{\textbf{u},\kappa^{\varepsilon}_{n}}_{s}\leq \max\limits_{k\in\mathcal{I}^{-\kappa^{\varepsilon}_{n}}}\left(Y^{\textbf{u},k}_s-\gamma_{\kappa^{\varepsilon}_{n},k}(s,X_{s})\right)+\frac{\varepsilon}{2^{n+1}}\right\}\wedge T\notag
\end{equation}
and
\begin{equation}
\kappa^{\varepsilon}_{n+1}=\arg\max_{k\in\mathcal{I}^{-\kappa^{\varepsilon}_{n}}}\left\{Y^{\textbf{u},k}_{\theta^{\varepsilon}_{n+1}}-\gamma_{\kappa^{\varepsilon}_{n},k}(\theta^{\varepsilon}_{n+1},X_{\theta^{\varepsilon}_{n+1}})\right\}.\notag
\end{equation}
In a similar manner as in the proof of Proposition \ref{admissible}, we can ensure that $\alpha^{\varepsilon}\in\mathcal{D}^{i}_{s}$ satisfies that $P[\theta^{\varepsilon}_{n}<T,\forall n \geq0]=0$. Let us prove now that $\expect\big[(\mathtt{C}^{\varepsilon}_{T})^{2}\big]<\infty$ and that $\alpha^{\varepsilon}$ is $\varepsilon$-optimal in $\mathcal{D}^{i}_{s}$ for the problem (\ref{5.11}). Following the strategy $\alpha^{\varepsilon}$ and since $(Y^{\textbf{u},i})_{i\in\mathcal{I}}$ solves the RBSDE (\ref{5.9}), it turns out that,
\begin{equation}\label{5.13}
Y^{\textbf{u},i}_{s}\leq Y^{\textbf{u},\kappa^{\varepsilon}_{n}}_{\theta^{\varepsilon}_n}+\displaystyle\int^{\theta^{\varepsilon}_{n}}_{s}f_{\alpha^{\varepsilon}}(r,X_{r},\textbf{u}_{r},Z^{\alpha^{\varepsilon}}_{r})dr-\mathtt{C}^{\alpha^{\varepsilon}}_{\theta^{\varepsilon}_{n}}-\int^{\theta^{\varepsilon}_{n}}_{s}Z^{\alpha^{\varepsilon}}_{r}dB_r+\varepsilon\sum^{n}_{i=1}\frac{1}{2^{i}},\hspace{0.2cm} \forall n\geq 1
\end{equation}
since $K^{\textbf{u},\kappa^{\varepsilon}_{n}}_{r}-K^{\textbf{u},\kappa^{\varepsilon}_{n}}_{\theta^{\varepsilon}_n}=0$ for $\theta^{\varepsilon}_{n}\leq r<\theta^{\varepsilon}_{n+1}$. Taking now the limit with respect to $n$ in (\ref{5.13}) we get:
\begin{equation}\label{5.14}
\begin{array}{ll}
Y^{\textbf{u},i}_{s} & \leq h_{\alpha^{\varepsilon}}(X_T)+\displaystyle\int^{T}_{s}f_{\alpha^{\varepsilon}}(r,X_{r},\textbf{u},Z^{\alpha^{\varepsilon}}_{r})dr-\mathtt{C}^{\alpha^{\varepsilon}}_{T}-\int^{T}_{s}Z^{\alpha^{\varepsilon}}_{r}dB_r+\varepsilon \\ \\
& = P^{\alpha^{\varepsilon}}_{s}-\mathtt{C}^{\alpha^{\varepsilon}}_{s}+\varepsilon.
\end{array} 
\end{equation}
Taking supremum over all $\alpha\in\mathcal{D}^{i}_{s}$, and next letting $\varepsilon\rightarrow 0$ and using the assumptions (\ref{H4}) and (\ref{H2})-(\ref{H2ii}),(\ref{H2iii}) satisfied for $h_i$ and $f_i$ respectively and since $\textbf{u}\in\mathcal{H}^{2,q}$, $Z^{\alpha^{\varepsilon}}\in\mathcal{H}^{2,d}$ and $(Y^{1}_{\cdot},\ldots,Y^{q}_{\cdot})\in(\mathscr{S}^{2})^{q}$, we deduce from (\ref{5.14}) that $\mathbb{E}[(\mathtt{C}^{\alpha^{\varepsilon}}_{T})^{2}]<\infty$. It follows that $Y^{\textbf{u},i}_{s}\leq \esssup_{\alpha\in\mathcal{D}^{i}_{s}}\left\{P^{\alpha}_{s}-\mathtt{C}^{\alpha}_{s}\right\}$. This last fact together with \eqref{41a} yield \eqref{5.11}. As a by-product, we obtain that the solution of (\ref{5.9}) is unique. $\hspace{1cm}\blacksquare$

\vspace{1cm}

Next for $\textbf{u}:=(u^{1},\ldots,u^{q})\in\mathcal{H}^{2,q}$ let us define
\begin{equation}
\Phi(\textbf{u}):= (Y^{\textbf{u},1}_{\cdot},\ldots,Y^{\textbf{u},q}_{\cdot}),\notag
\end{equation}
where $(Y^{\textbf{u},i}_{\cdot},Z^{\textbf{u},i}_{\cdot},K^{\textbf{u},i}_{\cdot})_{i=1,\ldots,q}$ is the solution of system (\ref{5.9}) which exists and is unique under the assumptions of Proposition \ref{5.3}. Note that the processes $(Y^{\textbf{u},i}_{\cdot},\ldots,Y^{\textbf{u},q}_{\cdot})$ belong to $(\mathscr{S}^{2})^{q}\subseteq\mathcal{H}^{2,q}$. Hence, $\Phi$ is a mapping from $\mathcal{H}^{2,q}$ to $\mathcal{H}^{2,q}$. 
\ms

We introduce the norm defined on $\mathcal{H}^{2,q}$ by
\begin{equation}
\left\|(u^{1},\ldots,u^{q})\right\|^{2}_{\beta}:=\expect\Bigg[\displaystyle\int^{T}_{0}e^{\beta s}\left(\sum^{q}_{i=1}\left|u^{i}_{s}\right|^{2}\right)ds\Bigg].\notag
\end{equation}
Note that $\left\|w\right\|_{\mathcal{H}^{2,q}}\leq\left\|w\right\|_{\beta}\leq e^{\beta T}\left\|w\right\|_{\mathcal{H}^{2,q}}$, for all $w\in\mathcal{H}^{2,q}$, implies that these norms are equivalent. For sake of completeness, we present the following result, established in Chassagneux et al. \cite{chassagneux2011note}, which ensures that $\Phi$ is a contraction on the Banach space $(\mathcal{H}^{2,q},\left\|\cdot\right\|_{\beta})$. 
\begin{PROP}\label{5.4}
Assume that for any $i,j\in\mathcal{I}$ the following hypotheses are in force:
\begin{enumerate}[label=(\roman*)]
\item $f_i$ verifies \textbf{\ref{H2}-\ref{H2ii},\ref{H2iii}};
\item $\gamma_{ij}$ (resp. $h_{i}$) verifies \textbf{\ref{H3}} (resp. \textbf{\ref{H4}}).
\end{enumerate}
Then, there exists $\beta_0\in\mathbb{R}$ such that the mapping $\Phi$ is a contraction operator on $(\mathcal{H}^{2,q},\left\|\cdot\right\|_{\beta_0})$. Therefore $\Phi$ has a fixed point $(Y^{1}_{\cdot},\ldots,Y^{q}_{\cdot})$ which belongs to $(\mathscr{S}^{2})^{q}$ and which provides a unique solution for system (\ref{RBSDEM5.4}).				
\end{PROP}
\textit{Proof.} Let $\textbf{u},\textbf{v}\in\mathcal{H}^{2,q}$ and consider the respective images under $\Phi$, $Y^{\textbf{u},i}_{\cdot}:=\Phi(\textbf{u})$ and $Y^{\textbf{v},i}_{\cdot}:=\Phi(\textbf{v})$. Besides, let us introduce the following ``auxiliary dominating'' RBSDE, for $i\in\mathcal{I}$:
\begin{equation}\label{mayorante}
\begin{cases}
\check{Y}^{i}_{s}=h_i(X_T)+\int^{T}_{s}\check{f}_i(r,X_{r},\check{Z}^{i}_r)dr+\check{K}^{i}_T-\check{K}^{i}_{s}-\int^{T}_{s}\check{Z}^{i}_rdB_r\hspace{0.5cm}\text{for all } 0\leq s\leq T, & \\ \check{Y}^{i}_{s}\geq \max\limits_{k\in\mathcal{I}^{-i}}\left\{\check{Y}^{k}_{s}-\gamma_{ik}(s,X_s)\right\}\hspace{0.5cm}\text{for all } 0\leq s\leq T, & \\
		\int^{T}_{0}\Big(\check{Y}^{i}_{s}-\max\limits_{k\in\mathcal{I}^{-i}}\left\{\check{Y}^{k}_{s}-\gamma_{ik}(s,X_{s})\right\}\Big)d\check{K}^{i,c}_{s}=0. & \\
		 \Delta_s\check{Y}^i_{\cdot}:=\check{Y}^{i}_s-\check{Y}^{i}_{s^{-}}=-\Big(\max\limits_{k\in\mathcal{I}^{-i}}\left\{\check{Y}^{k}_{s}-\gamma_{ik}(s,X_{s})\right\}-\check{Y}^{i}_s\Big)^{+}\hspace{0.5cm}\text{for all } 0\leq s\leq T 
\end{cases}
\end{equation}
where $\check{f}(s,X_{s},z^i)=\max\{f(s,X_{s},\textbf{u}_r,z^{i}),f(s,X_{s},\textbf{v}_r,z^{i})\}$, and $\check{K}^{i,c}_{\cdot}$ and $\check{K}^{i,d}_{\cdot}$ are the continuous and discontinuous parts of $\check{K}^{i}$. Note that by Proposition \ref{5.3} a unique solution $(\check{Y}^{i}_{\cdot},\check{Z}^{i}_{\cdot},\check{K}^{i}_{\cdot})$ exists for (\ref{mayorante}). 
\\
For $s\in[0,T]$ fixed, and for any $\alpha\in\mathcal{D}^{i}_{s}$, denote by $(U^{\alpha}_{\cdot},Z^{\alpha}_{\cdot})$,$(\bar{U}^{\alpha}_{\cdot},\bar{Z}^{\alpha}_{\cdot})$ and $(\check{U}^{\alpha}_{\cdot},\check{Z}^{\alpha}_{\cdot})$ the respective solutions of the following one-dimensional BSDEs: $\fr s\le T$,
\begin{equation}
\begin{array}{ll}
U^{\alpha}_{s}=h_{\alpha}(X_{T})+\int^{T}_{s}f_{\alpha}(r,X_{r},\textbf{u}_r,Z^{\alpha}_{r})dr-(\mathtt{C}^{\alpha}_{T}-\mathtt{C}^{\alpha}_{s})-\int^{T}_{s}Z^{\alpha}_{r}dB_r, & \\
\bar{U}^{\alpha}_{s}=h_{\alpha}(X_{T})+\int^{T}_{s}f_{\alpha}(r,X_{r},\textbf{v}_{r},\bar{Z}^{\alpha}_{r})dr-(\mathtt{C}^{\alpha}_{T}-\mathtt{C}^{\alpha}_{s})-\int^{T}_{s}\bar{Z}^{\alpha}_{r}dB_r,& \\
\check{U}^{\alpha}_{s}=h_{\alpha}(X_{T})+\int^{T}_{s}\check{f}_{\alpha}(r,X_{r},\check{Z}^{\alpha}_{r})dr-(\mathtt{C}^{\alpha}_{T}-\mathtt{C}^{\alpha}_{s})-\int^{T}_{s}\check{Z}^{\alpha}_{r}dB_r.
\end{array}\notag
\end{equation}
We deduce from Proposition \ref{5.3} that 
\begin{equation}\label{3.8}
Y^{\textbf{u},i}_{s}=\esssup_{\alpha\in\mathcal{D}^{i}_{s}}\left\{U^{\alpha}_{s}-\mathtt{C}^{\alpha}_{s}\right\},\hspace{0.2cm}\ Y^{\textbf{v},i}_{s}=\esssup_{\alpha\in\mathcal{D}^{i}_{s}}\left\{\bar{U}^{\alpha}_{s}-\mathtt{C}^{\alpha}_{s}\right\},\hspace{0.3cm} \check{Y}^{i}_{s}=\esssup_{\alpha\in\mathcal{D}^{i}_{s}}\left\{\check{U}^{\alpha}_{s}-\mathtt{C}^{\alpha}_{s}\right\}.
\end{equation}
Besides, note that for an $\varepsilon$-optimal strategy $\alpha^{\varepsilon}\in\mathcal{D}^{i}_{s}$, we have 
\begin{equation}\label{3.88}
 \check{Y}^{i}_{s}\leq \check{U}^{\alpha^\varepsilon}_s-\mathtt{C}^{\alpha^\varepsilon}_{s}+\varepsilon. 
\end{equation}
Using a comparison argument, we easily check that $\check{U}^{\alpha}_{\cdot}\geq U^{\alpha}_{\cdot}\vee \bar{U}^{\alpha}_{\cdot}$ for any strategy $\alpha\in\mathcal{D}^{i}_{s}$, and hence, by (\ref{3.8}) we get that $\check{Y}^{i}_{s}\geq Y^{\textbf{u},i}_{s}\vee Y^{\textbf{v},i}_{s}$. Therefore, taking into account the last two inequalities and (\ref{3.88}), we get that 
\begin{equation}
U^{\alpha^\varepsilon}_s-\mathtt{C}^{\alpha^{\varepsilon}}_{s}\leq Y^{\textbf{u},i}_s\leq \check{Y}^{i}_s\leq\check{U}^{\alpha^\varepsilon}_s-\mathtt{C}^{\alpha^{\varepsilon}}_{s}+\varepsilon\hspace{0.4cm} \text{and} \hspace{0.4cm}\bar{U}^{\alpha^\varepsilon}_s-\mathtt{C}^{\alpha^{\varepsilon}}_{s}\leq Y^{\textbf{v},i}_s\leq \check{Y}^{i}_s\leq\check{U}^{\alpha
^\varepsilon}_s-\mathtt{C}^{\alpha^{\varepsilon}}_{s}+\varepsilon.\notag
\end{equation}
This implies 
\begin{equation}
\left|Y^{\textbf{u},i}_s-Y^{\textbf{v},i}_s\right|\leq \left|\check{U}^{\alpha^\varepsilon}_s-U^{\alpha^\varepsilon}_s\right|+\left|\check{U}^{\alpha^\varepsilon}_s-\bar{U}^{\alpha^\varepsilon}_s\right|+2\varepsilon\notag,
\end{equation}
and by using the inequality $(a+b+c)^2\leq 4a^2+4b^2+2c^2$, we have
\begin{equation}\label{5.18}
\left|Y^{\textbf{u},i}_s-Y^{\textbf{v},i}_s\right|^{2}\leq 4\left|\check{U}^{\alpha^\varepsilon}_s-U^{\alpha^\varepsilon}_s\right|^{2}+4\left|\check{U}^{\alpha^\varepsilon}_s-\bar{U}^{\alpha^\varepsilon}_s\right|^{2}+4\varepsilon^{2}.
\end{equation}
Now, applying Ito's formula to $e^{\beta s}\left|\check{U}^{\alpha^\varepsilon}_s-U^{\alpha^\varepsilon}_s\right|^{2}$, using the inequality $\left|x\vee y-y\right|\leq\left|x-y\right|$ and the fact that $f_{\alpha_{\varepsilon}}$ is Lipschitz,  taking expectation, to obtain:  $\fr s\le T$,
\begin{equation}
\begin{array}{ll}
\expect\Big[e^{\beta s}\left|\check{U}^{\alpha^\varepsilon}_s-U^{\alpha^\varepsilon}_s\right|^{2}+\int^{T}_{s}e^{\beta r}\left|\check{Z}^{\alpha^\varepsilon}_r-Z^{\alpha^\varepsilon}_r\right|^{2}dr\Big] \leq -\expect\Big[\int^{T}_{s}\beta e^{\beta r}\left|\check{U}^{\alpha^\varepsilon}_r-U^{\alpha^\varepsilon}_r\right|^{2}dr\Big] \\ \\ \hspace{2cm}+2C\expect\Big[\int^{T}_{s}e^{\beta r}\left|\check{U}^{\alpha^\varepsilon}_r-U^{\alpha^\varepsilon}_r\right|\left(\left|\textbf{v}_r-\textbf{u}_r\right|+\left|\check{Z}^{\alpha^\varepsilon}_{r}-Z^{\alpha^\varepsilon}_{r}\right|\right)dr\Big].
\end{array}\notag
\end{equation}
The inequalities $2ab\leq \beta a^2+\frac{1}{\beta}b^2$ and $(a+b)^2\leq 2a^2 +2b^2$ also imply
\begin{equation}
\begin{array}{ll}
\expect\Big[e^{\beta s}\left|\check{U}^{\alpha^\varepsilon}_s-U^{\alpha^\varepsilon}_s\right|^{2}\Big]+\expect\Big[%
\int^{T}_{s}e^{\beta r}\left|\check{Z}^{\alpha^\varepsilon}_r-Z^{\alpha^\varepsilon}_r\right|^{2}dr\Big]\leq -\expect\Big[
\int^{T}_{s}\beta e^{\beta r}\left|\check{U}^{\alpha^\varepsilon}_r-U^{\alpha^\varepsilon}_r\right|^{2}dr\Big] \\ \\ \hspace{2cm} +\expect\Big[
\int^{T}_{s}\big\{\beta e^{\beta r}\left|\check{U}^{\alpha^\varepsilon}_r-U^{\alpha^\varepsilon}_r\right|^2+\frac{2C^2}{\beta}e^{\beta r}\left|\textbf{v}_r-\textbf{u}_r\right|^2+\frac{2C^2}{\beta}e^{\beta r}\left|\check{Z}^{\alpha^\varepsilon}_{r}-Z^{\alpha^\varepsilon}_{r}\right|^2\big\}dr\Big].\notag
\end{array}
\end{equation}
Rearranging terms, we obtain:
\begin{equation}
\expect\Big[e^{\beta s}\left|\check{U}^{\alpha^\varepsilon}_s-U^{\alpha^\varepsilon}_s\right|^{2}\Big]+\Big(1-\frac{2C^2}{\beta}\Big)\expect\bigg[
\int^{T}_{s} e^{\beta r}\left|\check{Z}^{\alpha^\varepsilon}_r-Z^{\alpha^\varepsilon}_r\right|^{2}dr\bigg]\leq \frac{2C^2}{\beta} \expect\Big[\int^{T}_{s} e^{\beta r}\left|\textbf{v}_r-\textbf{u}_r\right|^2dr\Big].\notag
\end{equation}
Taking $\beta\geq 2C^2$, we get 
\begin{equation}
\expect\Big[e^{\beta s}\left|\check{U}^{\alpha^\varepsilon}_s-U^{\alpha^\varepsilon}_s\right|^{2}\Big]\leq \frac{2C^2}{\beta}\expect\bigg[\int^{T}_{0} e^{\beta r}\left|\textbf{v}_r-\textbf{u}_r\right|^2dr\bigg].\notag
\end{equation}
Now, an analogous procedure to $e^{\beta s}\left|\check{U}^{\alpha_\varepsilon}_{s}-\bar{U}^{\alpha_\varepsilon}_{s}\right|^{2}$ lead to similar result, namely
\begin{equation}
\expect\Big[e^{\beta s}\left|\check{U}^{\alpha^\varepsilon}_s-\bar{U}^{\alpha^\varepsilon}_s\right|^{2}\Big]\leq \frac{2C^2}{\beta}\expect\bigg[\int^{T}_{0} e^{\beta r}\left|\textbf{v}_r-\textbf{u}_r\right|^2dr\bigg].\notag
\end{equation}
Combining these two inequalities with (\ref{5.18}), we deduce
\begin{equation}
\expect\Big[e^{\beta s}\left|Y^{\textbf{u},i}_s-Y^{\textbf{v},i}_s\right|^{2}\Big]\leq \frac{16C^2}{\beta}\expect\bigg[\int^{T}_{0} e^{\beta r}\left|\textbf{v}_r-\textbf{u}_r\right|^2dr\bigg]+4\varepsilon^2e^{\beta T}.\notag
\end{equation}
By integrating with respect to $s$ on both sides of the last inequality and taking into account the fact that such inequality holds true for any $i=1,\ldots,q$ and for all $s\in[0,T]$, we get 
\begin{equation}
\left\|\Phi(Y^{\textbf{u}})-\Phi(Y^{\textbf{v}})\right\|_{\beta}\leq4C\sqrt{Tq\beta^{-1}}\left\|\textbf{u}-\textbf{v}\right\|_{\beta}+2\varepsilon\sqrt{Tqe^{\beta T}}.\notag
\end{equation}
Finally, choosing $\beta_0>\max\left(16C^2Tq,2C^2\right)$ and taking $\varepsilon\rightarrow 0$, we see that this mapping is a contraction. This gives the existence and uniqueness of the system of RBSDE (\ref{RBSDEM5.4}). $\qed$

\section{The Markovian Framework}\label{sec4}
In this section we will provide more specifications to the process $X_{\cdot}$ treated in previous sections. Namely, we will assume now that this process has a Markovian evolution described by means of a stochastic differential equation (diffusion process) as in (\ref{eds}) below. Under this framework our previous analysis can be reduced to study a system of partial differential equations with obstacles (quasi-variational system). Among the main result of this section we can highlight the characterization of both the optimal function (\ref{pcg}) and the solution of the system of RBSDEs (\ref{RBSDEM5.4}) as a viscosity solution in a weak sense (see Theorem \ref{last}). We will start to introduce the following functions:
\begin{equation}
    \begin{cases}
      b:(t,x)\in\left[0,T\right]\times\mathbb{R}^{r}\mapsto b(t,x)\in\mathbb{R}^{r};& \\
        \sigma:(t,x)\in\left[0,T\right]\times\mathbb{R}^{r}\mapsto \sigma(t,x)\in\mathbb{R}^{r\times d},
  \end{cases}\notag
\end{equation}
satisfying the following hypotheses:
\ms

\noindent The functions $b$ and $\sigma$ are jointly continuous and Lipschitz continuous with respect to $x$ uniformly in $t$, that is, there exists a constant $C \geq 0$ such that for any $t \in[0,T ]$ and $x$, $x'\in\mathbb{R}^{r}$
    \begin{equation}\label{cdtlipschitz}
   \left|b(t,x)-b(t,x')\right|+\left|\sigma(t,x)-\sigma(t,x')\right|\leq C\left|x-x'\right|.
\end{equation}
Note that continuity and \eqref{cdtlipschitz} imply that $b$ and $\sigma$ are of linear growth, i.e., there exists a constant $C$ such that: \begin{equation}\label{cdtlinear}
 \left|b(t,x)\right|+\left|\sigma(t,x)\right|\leq C(1+\left|x\right|),\,\,\forall (t,x) .
\end{equation}
It is well known that under \eqref{cdtlipschitz}-\eqref{cdtlinear},
there exists a unique Markov process $(X^{t,x}_{s})_{s\leq T}$, for
$\left(t,x\right)\in\left[0,T\right]\times\mathbb{R}^{r}$, that is a
(strong) solution of the following standard stochastic differential
equation:
\begin{equation}\label{eds}
    \begin{cases}
dX^{t,x}_{s}=b(s,X^{t,x}_{s})ds+\sigma(s,X^{t,x}_{s})dB_s\hspace{0.5cm}\text{for all } t\leq s\leq T;& \\ X^{t,x}_{s}=x \hspace{0.5cm}\text{for all } 0\leq s\leq t,
\end{cases}
\end{equation}
satisfying the following estimates: For any $p\geq 2$, $x,x'\in \mathbb{R}^r$ and $s\ge t$
    \begin{align}\label{5.33}
& \expect\big[\sup_{s\leq T}\left|X^{t,x}_{s}\right|^{p}\big]\leq C\left(1+\left|x\right|^{p}\right),\quad \expect [\sup_{r\in [t,s]}|X^{t,x}_r-x|^p]\leq M_p(s-t)(1+|x|^p)] \mbox{ and } \\&\quad\qquad
\expect [\sup_{r\in [t,s]}|X^{t,x}_r-X^{t,x'}_r-(x-x')|^p]\leq M_p(s-t)|x-x'|^p \nonumber\end{align}for some constant $M_p$ (one can see Karatzas and Shreve \cite{karashreve} or Revuz and Yor \cite{revuzyor}, for more details).

Recall that the associated infinitesimal generator to $(X^{t,x}_{s})_{s\leq T}$ is given by :
\begin{equation}
    \mathcal{L}\phi(t,x)=\frac{1}{2}Tr\left[\left(\sigma.\sigma^{\mathtt{T}}\right)(t,x)D^{2}_{xx}\phi(t,x)\right]+b(t,x)^{\mathtt{T}}D_x\phi(t,x)\notag
    \end{equation}
for $\phi$ in $C^{1,2}([0,T]\times\mathbb{R}^{r};\mathbb{R})$ ($Tr(.)$ is the trace of a square matrix and, $\mathtt{A}^{\mathtt{T}}$ is the transpose of a matrix $\mathtt{A}$).
\bigskip

Now let $(t,x)\in \left[0,T\right]\times\mathbb{R}^{r}$ be fixed and let $((Y^{t,x,i}_{s},Z^{t,x,i}_{s},K^{t,x,i}_{s})_{t\leq s\leq T})_{i=1,\ldots,q}$ be the unique solution of system \eqref{RBSDEM5.4} when the process $X$ is taken to be equal to $X^{t,x}$ of \eqref{eds}, i.e., the solution associated with
$(f_i(s,X^{t,x}_{s},y^1,\ldots,y^i,\ldots y^q,z^i),h_i(X^{t,x}_{T}),g_{ik}(s,X^{t,x}_{s}))$ ($g_{ik}$ are the switching costs and they satisfy the same assumptions as $\gamma_{ik}$ in Assumption (\ref{assH})) with $y^{i}\in\mathbb{R}$ and $z^i\in\mathbb{R}^d$.

Assume now that Assumptions (\ref{assH}) are satisfied. Since we are in the Markovian framework then there exist deterministic functions $u^i$, $\ici$, with polynomial growth such that for any $(t,x)$
$$
Y^{t,x,i}_s=u^i(s,X^{t,x}_s),\,\,\ici, \ \p-a.s., \,\,\fr s\in [t,T]. $$
Note that we also have \begin{equation}\label{deterministic}
u^i(t,x)= Y^{t,x,i}_{t}, \hspace{0.2cm}\text{for}\hspace{0.2cm} (t,x)\in[0,T]\times\mathbb{R}^{r}\hspace{0.2cm}\text{and}\hspace{0.2cm} \ici .
\end{equation}
On the other hand the polynomial growth of $u^i$ stems from the polynomial growths of the data assumed in Assumption (\ref{assH}) and the BSDEs \eqref{BSDEMS}, \eqref{BSDEMM} as well.

\emph{Notation:} For a sake of simplicity of notation, hereafter we sometimes denote by $(\psi)_{k=1,\dots, q}:= (\psi_1,\dots,\psi_q)$, for some generic function or vector $\psi$.

\begin{REM}\label{fimonotone}
From now on we will assume that $f_i$ is non-decreasing w.r.t $y^k$ for any $k=1,...,q$ and not only w.r.t $y^1, ...,y^{i-1},y^{i+1},...,y^q$ (as precised in  \ref{H2}-\ref{H2iv}). This assumption is not really restrictive since by considering the system of RBSDEs  verified by $(e^{\alpha t}Y^i_t)_{t\le T}$ , we obtain new generators $F_i$ given by 
$$
F_i(t,y^1, ...,y^m,z^i):=e^{\alpha t}f_i(t,x,e^{-\alpha t}y^1,...,e^{-\alpha t}y^m,e^{-\alpha t}z^i)-\alpha y^i
$$
which have the same properties as $(f_i)_{i=1,\dots,q}$. Moreover, with an appropriate choice of $\alpha$, those new generators are non-decreasing w.r.t $y^k$ for any $k=1,\dots,q$, i.e., they fulfill the property we are requiring for $(f_i)_{i=1,\dots,q}$ (one can see Hamad\`ene and Morlais \cite{hamadene2013viscosity}, for more details on this transform). 
\end{REM}
\ms

Our main interest will be to show that the function
$(u^i)_{i=1,\dots,q}:(t,x)\in\left[0,T\right]\times\mathbb{R}^{r}\mapsto(u^i(t,x))_{i=1,\dots,q}\in\mathbb{R}^{q}$
is a solution in a weak viscosity sense for the Hamilton-Jacobi-Bellman
system of PDEs associated with the switching problem. In the case
when the functions $g_{ij}$ and $h_i$, $i,j\in \mathcal{I}$, are
continuous, this system reads as: for all $i\in\mathcal{I}$,
\begin{equation}\label{pdes4.2}
    \begin{cases}
      \min \{v^i(t,x)-\max\limits_{k\in\mathcal{I}^{-i}}\left(v^{k}(t,x)-g_{ik}(t,x)\right);
      \hspace{0.1cm} -\partial_tv^i(t,x)-\mathcal{L}v^i(t,x)-& \\ \hspace{2cm}-f_i(t,x,\left(v^{1},
      \ldots,v^{i},\ldots,v^{q}\right)(t,x),\sigma^{\mathtt{T}}(t,x)D_xv^{i}(t,x))\}=0; & \\
      v^i(T,x)=h_i(x).
         \end{cases}
\end{equation}and it is shown that $(u^i)_{i=1,\dots,q}$ is the unique viscosity solution of system \eqref{pdes4.2}. But in our framework the functions $g_{ij}$, $i,j\in \mathcal{I}$,
are no longer continuous w.r.t $t$, therefore the definition should
be adapted. We are going to show that $(u^i)_{i=1,\dots,q}$ is a viscosity
solution in a weak sense for the HJB system of PDEs \eqref{pdes4.2}, associated with
the swiching problem, and which we are going to define in what
follows. This definition is inspired by Ishii's works 
\cite{{ishii1987perron},{ishii1987Hamilton}}, and also by the paper of Barles and
Perthame \cite{Barles1987Discontinuous}.

To proceed for a locally bounded $\mathbb{R}$-valued function
$v(t,x)$, $(t,x)\in [0,T]\times \mathbb{R}^\ell$ ($\ell \ge 1$), we
define its lower (resp. upper) semi-continuous envelope $v_*$ (resp.
$v^*$) as follows: For any $(t,x)\in [0,T]\times \mathbb{R}^\ell$,
\begin{equation}
v^{\ast}(t,x):= \limsup_{(t',x')\rightarrow
(t,x),\hspace{0.1cm}t'<T}v(t',x') \hspace{0.3cm}\text{(resp.
}\hspace{0.3cm} v_{\ast}(t,x):= \liminf_{(t',x')\rightarrow
(t,x),\hspace{0.1cm}t'<T}v(t',x')).\notag
\end{equation}
Note that the function $v^*$ (resp. $v_*$) can also be seen as the  smallest usc (resp. lsc) function which is greater (resp. smaller) than $v$. On the other hand, the following properties of the semi-continuous
envelopes of functions will be useful later. 
\begin{LEM}\label{semicontinuity} Let $(t,x)\in
[0,T]\times \mathbb{R}^\ell$ and $\varphi_i(t,x)$, $i=1,2$, be two
locally bounded $\mathbb{R}$-valued functions. We then have:
\begin{enumerate}[label=(\roman*)]
\item If $\varphi_1$ is continuous then $(\varphi_1+\varphi_2)_*=
\varphi_1+(\varphi_2)_*$ and $(\varphi_1+\varphi_2)^*=
\varphi_1+(\varphi_2)^*$.

\item $(-\varphi_1)_*=-(\varphi_1)^*$. 

\item $(\varphi_1 \wedge \varphi_2)_*=(\varphi_1)_* \wedge (\varphi_2)_*$ and 
$(\varphi_1 \vee \varphi_2)^*=(\varphi_1)^* \vee (\varphi_2)^*$.

\item If 
$\varphi_1$ is continuous then $(\varphi_1 \wedge \varphi_2)^*=\varphi_1\wedge (\varphi_2)^*$ and 
$(\varphi_1 \vee \varphi_2)_*=\varphi_1\vee (\varphi_2)_*$.
\end{enumerate}
\end{LEM}
\noindent \emph{Proof.} (i) Obviously we have $\varphi_1+\varphi_2\ge \varphi_1+(\varphi_2)_*$ and then 
$(\varphi_1+\varphi_2)_*\ge \varphi_1+(\varphi_2)_*$ since this latter is lsc. On the other hand 
$(\varphi_1+\varphi_2)_*-\varphi_1\le \varphi_2$ and then $(\varphi_1+\varphi_2)_*-\varphi_1\le (\varphi_2)_*$ since $(\varphi_1+\varphi_2)_*-\varphi_1$ is lsc. This completes the proof of the claim as the other property can be obtained similarly. 

Points (ii) and (iii) are rather obvious, we then leave their proofs to the care of the reader.

(iv) First note that $(\varphi_1 \wedge \varphi_2)^*\le \varphi_1\wedge (\varphi_2)^*$. Next let $((t_n,x_n))_n$ be a sequence such that $(\varphi_2(t_n,x_n))_n\rightarrow (\varphi_2)^*(t,x)$ as $n\to\infty$. As $\varphi_1$ is continuous then 
$\varphi_1(t_n,x_n)\wedge \varphi_2(t_n,x_n)\rightarrow \varphi_1(t,x)\wedge (\varphi_2)^*(t,x)$ as $n\to\infty$. Therefore, by definition of the usc envelope, $(\varphi_1 \wedge \varphi_2)^*(t,x)\ge \varphi_1(t,x)\wedge (\varphi_2)^*(t,x)$ which completes the proof of the first claim. The proof of the other one is similar.  $\qed$
\bigskip

Next for $i=1,...,q$, let us denote by $F_i$ the non-linearity which
defines the $i-th$ equation in \eqref{pdes4.2}, i.e.,
\begin{equation}
F_{i}(t,x,(y^{j})_{j=1,\ldots,q},r,p,X)=\min\big\{y^{i}-\max\limits_{k\in\mathcal{I}^{-i}}(y^{k}-g_{ik}(t,x));\hspace{0.1cm}G_i(t,x,(y^{j})_{j=1,\ldots,q},r,p,X)
\big\}
\end{equation}where
\begin{equation}
G_i(t,x,(y^{j})_{j=1,\ldots,q},r_{i},p_{i},X_{i})=-r_{i}-\frac{1}{2}
Tr(\sigma\sigma^{\mathtt{T}}X_{i})-b^{\mathtt{T}}p_{i}-f_i(t,x,(y^{j})_{j=1,\ldots,q},\sigma^{\mathtt{T}}p_i).
\end{equation}
Note that by Assumption (\ref{assH}) on $f_i$ and \eqref{cdtlipschitz}, the
function $G_i$ is jointly continuous in its arguments. Therefore, taking into account the results of Lemma \ref{semicontinuity}, for any $i=1,\dots,q$, the semi-continuous envelopes of
$F_i$ (in all arguments) are given by:
$$
F_i^*(t,x,(y^{j})_{j=1,\ldots,q},r,p,X)=
\min\big\{y^{i}-(\max\limits_{k\in\mathcal{I}^{-i}}(y^{k}-g_{ik}(t,x)))_*;\hspace{0.1cm}G_i(t,x,(y^{j})_{j=1,\ldots,q},r,p,X)
\big\}
$$ and
$$
(F_i)_*(t,x,(y^{j})_{j=1,\ldots,q},r,p,X)=
\min\big\{y^{i}-(\max\limits_{k\in\mathcal{I}^{-i}}(y^{k}-g_{ik}(t,x)))^*;\hspace{0.1cm}G_i(t,x,(y^{j})_{j=1,\ldots,q},r,p,X)
\big\}.
$$

We are now ready to precise the definition of the viscosity solution of HJB system associated with the switching problem. As noticed previously it is inspired by the papers \cite{Barles1987Discontinuous,{ishii1987perron},{ishii1987Hamilton}}. On the other hand, the discontinuities of the functions $(u^i)_{i=1,\dots,q}$ generated by the ones of $(g_{ij})_{i,j\in \mathcal{I}}$ make that the terminal condition 
at time $t=T$ is not the same as in \eqref{pdes4.2}, but should be adapted as well to this weak sense (see e.g. \cite{Bouchard2009Stochastic}). 

\begin{DEF}\label{definitionvs}
 Let $\mathbf{v}:=(v^1,\ldots,v^q)$ be  a locally bounded function from $[0,T]\times\mathbb{R}^{r}$ into $\mathbb{R}^{q}$. 
 \begin{description}
\item[(1)] We say that $\mathbf{v}$ is a viscosity subsolution of \eqref{pdes4.2} if for any $i\in\mathcal{I}$, and $x_0\in\mathbb{R}$,
 \begin{description}
\item[(a)]  $v^{i*}$ verifies the following inequality at point $(T,x_0)$:
\begin{equation}\label{defsoussolviscoenT}
\min\big\{v^{i\ast}(T,x_0)-h_{i}(x_0);\hspace{0.2cm} u^{i\ast}(T,x_0)-\max\limits_{j\in\mathcal{I}^{-i}}\left(v^{j\ast}-g_{ij}\right)^{\ast}(T,x_0)\big\}\leq 0.
\end{equation}
\item[(b)] Moreover, at $(t_0,x_0)\in [0,T)\times \mathbb{R}^{r}$, the function $v^i$ is such that,
for and any $\phi\in C^{1,2}([0,T]\times\mathbb{R}^{r})$ with $\phi(t_0,x_0)=v^{i\ast}(t_0,x_0)$ and $\phi-v^{i\ast}$ attaining its minimum at $(t_0,x_0)$, we have
\begin{equation}\begin{array}{l}
 (F_i)_*(t_0,x_0,(v^{j\ast}(t_0,x_0))_{j=1,\ldots,q},\partial_t\phi(t_0,x_0),D_x\phi(t_0,x_0),D^2_{xx}\phi(t_0,x_0))\\\qquad \qquad = \min \bigg\{v^{i\ast}(t_0,x_0)-\max\limits_{k\in\mathcal{I}^{-i}}\left(v^{k*}-g_{ik}\right)^{\ast}(t_0,x_0);\\\qquad\qquad -\left(\partial_t+\mathcal{L}\right)\phi(t_0,x_0) -f_i(t_0,x_0,v^{j*}(t_0,x_0))_{j=1,\ldots,q}(t_0,x_0),\left(\sigma^{\mathtt{T}}D_x\right)\phi(t_0,x_0))\bigg\}\leq 0.\end{array}\notag
\end{equation}
\end{description}
\item[(2)] In the same manner, $\mathbf{v}$ is said to be a viscosity supersolution of \eqref{pdes4.2} if for any $i\in\mathcal{I}$, and $x_0\in\mathbb{R}$,
 \begin{description}
\item[(a)] $v^i_*$ verifies at $(T,x_0)$ the following:
\begin{equation}\label{defsursolviscoenT}
\min\Big\{v^{i}_{\ast}(T,x_0)-h_{i}(x_0);\hspace{0.2cm} v^{i}_{\ast}(T,x_0)-\big(\max\limits_{j\in\mathcal{I}^{-i}}\left(v^{j}_{\ast}-g_{ij}\right)\big)_{\ast}(T,x_0)\Big\}\geq 0.
\end{equation}
\item[(b)] Similarly, at $(t_0,x_0)\in [0,T)\times \mathbb{R}^{r}$, $v^i$ satisfies the next: for and any $\phi\in C^{1,2}([0,T]\times\mathbb{R}^{r})$ with $\phi(t_0,x_0)=v^{i}_*(t_0,x_0)$ and $\phi-v^{i}_*$ attaining its maximum at $(t_0,x_0)$, we have
\begin{equation}\begin{array}{l}
 (F_i)^*(t_0,x_0,(v^{j}_*(t_0,x_0))_{j=1,\ldots,q},\partial_t\phi(t_0,x_0),D_x\phi(t_0,x_0),D^2_{xx}\phi(t_0,x_0))\\\qquad = \min \bigg\{v^{i}_*(t_0,x_0)-(\max\limits_{k\in\mathcal{I}^{-i}}\left(v^{k}_*-g_{ik}\right))_{\ast}(t_0,x_0);\\\qquad\qquad -\left(\partial_t+\mathcal{L}\right)\phi(t_0,x_0) -f_i(t_0,x_0,(v^{j}_*(t_0,x_0))_{j=1,\dots,q},\left(\sigma^{\mathtt{T}}D_x\right)\phi(t_0,x_0))\bigg\}\geq 0.\end{array}\notag
\end{equation}
\end{description}
\item[(3)] We say that $\mathbf{v}$ is viscosity solution of (\ref{pdes4.2}) if it is both  a viscosity sub. and supersolution. 
\end{description}
\end{DEF}
\noindent

To proceed we are going to show that the functions $(u^i)_{i=1,...,q}$ is a viscosity solution of the system 
\eqref{pdes4.2} in a weak sense, i.e., according to Definition \ref{definitionvs}. However we need some preliminary results which we give as lemmas hereafter. From now $B_{\eta}(t_0,x_0)$ is the open ball of radius $\eta$ and center $(t_0,x_0)$.
\begin{LEM}\label{lemaalex}
Under the Assumption (\ref{H2}), the mapping
$$
(t,x)\longmapsto f_i\big(t,x,\big(v^{1\ast},\ldots,v^{q\ast}\big)(t,x),\left(\sigma^{\mathtt{T}}D_x\right)\phi(t,x)\big)
$$
is u.s.c. for any $\phi \in C^{1,2}([0,T]\times\mathbb{R}^r)$.
\end{LEM}
\noindent \textit{Proof.} Let $(t_0,x_0)\in [0,T]\times\mathbb{R}^{r}$. Since $v^{k\ast}$ is \textit{u.s.c} for $k=1,\dots,q$, then for all $\varepsilon>0$  there exists $\eta_{\varepsilon}>0$ such that for all $(t,x)$, satisfying $\left\|(t,x)-(t_0,x_0)\right\|<\eta_{\varepsilon}$, we have
\begin{equation}
v^{k\ast}(t_0,x_0)\geq v^{k\ast}(t,x)-\varepsilon \hspace{0.7cm}\text{for all}\hspace{0.3cm} k=1,\ldots,q.\notag
\end{equation}
Next, by monotonicity and Lipschitz properties of $f_i$, 
 for all $(t,x)\in B_{\eta_{\varepsilon}}(t_0,x_0)$ we get
\begin{equation}
\begin{split}
&\ f_i\big(t_0,x_0,\big(v^{k\ast}\big)_{k=1,\ldots,q}(t_0,x_0),\left(\sigma^{\mathtt{T}}D_x\right)\phi(t_0,x_0)\big) \geq f_i\big(t_0,x_0,\big(v^{k\ast}(t,x)-\varepsilon\big)_{k=1,\ldots,q},\big(\sigma^{\mathtt{T}}D_x\big)\phi(t_0,x_0)\big)\\ 
&\ \geq f_i\big(t_0,x_0,\big(v^{k\ast}(t,x)\big)_{k=1,\dots,q},\big(\sigma^{\mathtt{T}}D_x\big)\phi(t_0,x_0)\big)-C\varepsilon =f_i(t,x,\big(v^{k\ast}(t,x)\big)_{k=1,\ldots,q},\left(\sigma^{\mathtt{T}}D_x\right)\phi(t,x))-C\varepsilon+\\  &\ +\Big\{f_i\big(t_0,x_0,\big(v^{k\ast}(t,x)\big)_{k=1,\ldots,q},\big(\sigma^{\mathtt{T}}D_x\big)\phi(t_0,x_0)\big)  -f_i(t,x,\big(v^{k\ast}(t,x)\big)_{k=1,\ldots,q},\big(\sigma^{\mathtt{T}}D_x\big)\phi(t,x))\Big\},
\end{split}\notag
\end{equation}
where $C$ is the Lipschitz constant of $f_i$. By continuity of $f_i$ with respect to $(t,x)$ and Lipschitz in $z^i$ (Assumptions (\ref{H2})-(\ref{H2i}) and (\ref{H2})-(\ref{H2ii})), the quantity inside the brackets goes to zero as $(t,x)\rightarrow (t_0,x_0)$. Therefore, taking a suitable $\eta_{\varepsilon}>0$, we obtain
\begin{equation}
f_i\big(t_0,x_0,\big(v^{k\ast}(t_0,x_0)\big)_{k=1,\ldots,q},\left(\sigma^{\mathtt{T}}D_x\right)\phi(t_0,x_0)\big)\geq f_i(t,x,\big(v^{k\ast}(t,x)\big)_{k=1,\ldots,q},\left(\sigma^{\mathtt{T}}D_x\right)\phi(t,x))-C'\varepsilon\notag
\end{equation}
for all $(t,x)\in B_{\eta_{\varepsilon}}(t_0,x_0)$ and for some other constant $C'$ and the claim follows.\hspace{3cm}$\blacksquare$
\begin{LEM}\label{lemasaid} Let $\phi\in C^{1,2}([0,T]\times\mathbb{R}^{r})$, $(t_0,x_0)\in [0,T)\times\mathbb{R}^{r}$ and $\phi(t_0,x_0)=v^{i\ast}(t_0,x_0)$. If
\begin{equation}\label{72a}
\phi(t_0,x_0)=v^{i\ast}(t_0,x_0)>\max\limits_{k\in\mathcal{I}^{-i}}\big(v^{k\ast}-g_{ik}\big)^{\ast}(t_0,x_0)
\end{equation}
and
\begin{equation}
-\left(\partial_t+\mathcal{L}\right)\phi(t_0,x_0)>f_{i}(t_0,x_0,(v^{k\ast})_{k=1,\dots,q}(t_0,x_0),(\sigma^{\mathtt{T}}D_x)\phi(t_0,x_0)),
\end{equation}
then there exist $\varepsilon$ and a ball $B_{\eta_{\varepsilon}}(t_0,x_0)$ such that for all $(t,x)\in B_{\eta_{\varepsilon}}(t_0,x_0)$ we have:
\begin{equation}\label{4}
\phi(t,x)\geq\max\limits_{k\in\mathcal{I}^{-i}}\big(v^{k\ast}(t,x)-g_{ik}(t,x)\big)+\varepsilon
\end{equation}
and
\begin{equation}\label{5}
-\left(\partial_t+\mathcal{L}\right)\phi(t,x)\geq f_{i}\big(t,x,(v^{k\ast})_{k=1,\dots,q}(t,x),(\sigma^{\mathtt{T}}D_x)\phi(t,x)\big)+\varepsilon.
\end{equation}
\end{LEM}
\noindent \textit{Proof}: By (\ref{72a}) and the continuity of $\phi$ there exist $\varepsilon$ and a ball $B_{\eta_{\varepsilon}}(t_0,x_0)$ such that
\begin{equation}\label{1}
\phi(t,x)\geq\max\limits_{k\in\mathcal{I}^{-i}}\big(v^{k\ast}(t_0,x_0)-g_{ik}(t_0,x_0)\big)^{\ast}+2\varepsilon
\end{equation}
for all $(t,x)\in B_{\eta_{\varepsilon}}(t_0,x_0)$. Next, by the \textit{u.s.c} property, there exists $\eta_{\varepsilon}^{'}$ such that for all $(t,x)\in B_{\eta^{'}_{\varepsilon}}(t_0,x_0)$ we have
\begin{equation}\label{2}
\begin{array}{ll}
\max\limits_{k\in\mathcal{I}^{-i}}\left(v^{k\ast}(t_0,x_0)-g_{ik}(t_0,x_0)\right)^{\ast}\geq\max\limits_{k\in\mathcal{I}^{-i}}\left(v^{k\ast}(t,x)-g_{ik}(t,x)\right)^{\ast}-\varepsilon & \\ & \\ \hspace{5.4cm} \geq \max\limits_{k\in\mathcal{I}^{-i}}\left(v^{k\ast}(t,x)-g_{ik}(t,x)\right)-\varepsilon
\end{array}
\end{equation}
where in the last inequality we use that the usc envelope of a function is greater or equal to the function itself. Therefore, from (\ref{1}), (\ref{2}) and assuming , without loss of generality, that $\eta_{\varepsilon}\leq \eta_{\varepsilon}^{'}$  we have
\begin{equation}
\phi(t,x)\geq\max\limits_{k\in\mathcal{I}^{-i}}\big(v^{k\ast}(t,x)-g_{ik}(t,x)\big)+\varepsilon
\end{equation}
for all $(t,x)\in B_{\eta_{\varepsilon}}(t_0,x_0)$. 

As for the second inequality we can do a similar procedure since $\left(\partial_t+\mathcal{L}\right)\phi$ is continuous and $(t,x)\mapsto f_{i}(t,x,(v^{k\ast})_{k=1,\dots,q}(t,x),(\sigma^{\mathtt{T}}D_x)\phi(t,x))$ is \textit{u.s.c}. Namely, there exist $\varepsilon^{'}$ and $\eta_{\varepsilon}^{''}$ such that for each $(t,x)\in B_{\eta_{\varepsilon}^{''}}(t_0,x_0)$ we have
\begin{equation}
-\left(\partial_t+\mathcal{L}\right)\phi(t,x)\geq f_{i}\big(t,x,(v^{k\ast})_{k=1,\dots,q}(t,x),(\sigma^{\mathtt{T}}D_x)\phi(t,x)\big)+\varepsilon'.
\end{equation}
Now, supposing, without loss of generality, that $\varepsilon\leq\varepsilon'$ and $\eta_{\varepsilon}\leq\eta_{\varepsilon}^{''}$, we have that inequalities (\ref{4}) and (\ref{5}) hold true for all $(t,x)\in B_{\eta_{\varepsilon}}$.
\begin{REM}\label{remark55}
In a similar manner, it is possible to obtain a parallel result as in Lemmas \ref{lemaalex} and \ref{lemasaid} for $v^{i}_{\ast}$ in lieu of $v^{i\ast}$. Namely, it can be proved that under Assumption (\ref{H2}) the mapping 
$$
(t,x)\longmapsto f_i\big(t,x,(v^{k}_*)_{k=1,\dots,q}(t,x),\big(\sigma^{\mathtt{T}}D_x\big)\phi(t,x)\big)
$$
is l.s.c., and if
\begin{equation}
-\left(\partial_t+\mathcal{L}\right)\phi(t_0,x_0)<f_{i}(t_0,x_0,(v^k_{\ast})_{k=1,\dots,q}(t_0,x_0),(\sigma^{\mathtt{T}}D_x)\phi(t_0,x_0)),
\end{equation}
then there exists $\varepsilon>0$ and $\eta_{\varepsilon}$ such that for all $(t,x)\in B_{\eta_{\varepsilon}}(t_0,x_0)$:
\begin{equation}
-\left(\partial_t+\mathcal{L}\right)\phi(t,x)\leq f_i\big(t,x,(v^{k}_*)_{k=1,\dots,q}(t,x),\big(\sigma^{\mathtt{T}}D_x\big)\phi(t,x))-\varepsilon.\notag
\end{equation}
The proofs are very similar as the proofs given in the aforementioned lemmas, so shall omit them.
\end{REM}
Finally, we recall two comparison results for BSDE and RBSDE that we have borrowed from Lemma 4.1 and Proposition 4.2, in Dumitrescu et al.  \cite{dumitrescu2016weak}.
\begin{LEM}\label{lemma46} Fix $t_0\in[0,T]$ and let $\theta$ be a stopping time with values in $[t_0,T]$. Consider two random variables $\xi_1$ and $\xi_2\in L^{2}(\mathcal{F}_\theta)$ and two  drivers (a.k.a generators) $f_1$, $f_2$ such that $f_2$ satisfies \ref{H2} with Lipschitz constant $C>0$. For $i=1,2$, let $(Y^{i}_{t},Z^{i}_{t})$ be the solution in $\mathscr{S}^{2}\times \mathcal{H}^2$ of the BSDE with associated data  $(f_i,\xi_i)$, and terminal time $\theta$. In this case, $f_i$ and $\xi_i$ represent the driver and terminal condition, respectively. Suppose that for some $\epsilon>0$ we have
\begin{equation}
1_{\{t_0\leq t\leq \theta\}}(t)f_1(t,Y^{1}_{t},Z^{1}_{t})\geq 1_{\{t_0\leq t\leq \theta\}}(t)f_2(t,Y^{1}_{t},Z^{1}_{t}),\,\,dt\otimes d\P-a.e.\text{ and } \xi_1\geq\xi_2+\varepsilon, \,\,\P- a.s.\notag
\end{equation}
Then we have $Y^{1}_{t}\geq Y^{2}_{t}+\varepsilon e^{-CT}$, $\P$-a.s. for each $t\in[t_0,\theta]$.
\end{LEM}
\begin{LEM}[A comparison result between a BSDE and a RBSDE]\label{compareBR}
 Fix $t_0\in[0,T]$ and let $\theta$ be a stopping time on $[t_0,T]$. Consider the random variable $\xi_1\in L^2(\mathcal{F}_\theta)$  and a driver $f_1$. Let $(Y^{1}_{t},Z^{1}_{t})$ be the associated BSDE solution  with driver $f_1$, terminal time $\theta$ and terminal condition $\xi_1$. Consider also $g_{2}(\cdot)\in\mathscr{S}^2$ and let $f_2$ be a driver satisfying \ref{H2} with Lipschitz constant $C>0$. Assume the existence of the solution $Y^{2}_{t}$ of the associated RBSDE with driver $f_2$, terminal time $\theta$ and obstacle $g_{2}$, and assume that
\begin{equation}\begin{array}{l}
1_{\{t_0\leq t\leq \theta\}}(t)f_1(t,Y^{1}_{t},Z^{1}_{t})\geq 1_{\{t_0\leq t\leq \theta\}}(t)f_2(t,Y^{1}_{t},Z^{1}_{t}),\,\,dt\otimes d\P-a.e.\\\text{ and } \\ 1_{\{t_0\leq t\leq \theta\}}(t)Y^{1}_{t}\geq 1_{\{t_0\leq t\leq \theta\}}(t)(g_{2}({t})+\varepsilon), \fr t\ge 0, \P-a.s.\notag\end{array}
\end{equation}
where $\varepsilon$ is a positive constant. Then,  we have $Y^{1}_{t}\geq Y^{2}_{t}+\varepsilon e^{-CT}$, $\P-a.s.$, for each $t\in[t_0,\theta]$.
\end{LEM}
We now give the main result of this section.
\begin{THM}\label{last}
The function $\mathbf{u}:=(u^{1},\ldots,u^{q})$, where for each $i=1,\ldots,q$, $u^i$ is defined as in (\ref{deterministic}), is a weak viscosity solution of the system (\ref{pdes4.2}).
\end{THM}
\noindent \textit{Proof.} \underline{Step 1}: Viscosity sub-solution property on $[0,T)\times \mathbb{R}^r$. 
\bs 

\noindent Let $\phi\in C^{1,2}([0,T]\times\mathbb{R}^{r})$ and $(t_0,x_0)\in[0,T)\times\mathbb{R}^{r}$ be such that $\phi(t,x)\geq u^{i\ast}(t,x)$, for all $(t,x)\in[0,T)\times\mathbb{R}^{r}$ and $\phi(t_0,x_0)=u^{i\ast}(t_0,x_0$). Without loss of generality, we can assume that the minimum of $\phi-u^{i\ast}$ attained at $(t_0,x_0)$ is strict. We need to show that if
\begin{equation}
\phi(t_0,x_0)=u^{i\ast}(t_0,x_0)>\max\limits_{k\in\mathcal{I}^{-i}}\big(u^{k\ast}(t_0,x_0)-g_{ik}(t_0,x_0)\big)^{\ast}
\end{equation}
 then
\begin{equation}\label{subsolutionproperty}
 -\left(\partial_t+\mathcal{L}\right)\phi(t_0,x_0)-f_i\big(t_0,x_0,(u^{k*})_{k=1,\dots,q}(t_0,x_0),(\sigma^{\mathtt{T}}D_x)\phi(t_0,x_0)\big)\leq 0.
\end{equation}
We proceed by contradiction; i.e. we shall assume \begin{equation}
 -\left(\partial_t+\mathcal{L}\right)\phi(t_0,x_0)-f_i\big(t_0,x_0,(u^{k*})_{k=1,\dots,q}(t_0,x_0),(\sigma^{\mathtt{T}}D_x)\phi(t_0,x_0)\big)>0, \notag
\end{equation}
then by Lemma \ref{lemasaid} there exists $\varepsilon>0$ and $\eta_\varepsilon>0$ such that, for all $(t,x)\in B_{\eta_{\varepsilon}}(t_0,x_0)$, we have both
\begin{equation}\label{weak33}
\phi(t,x)\geq\max\limits_{k\in\mathcal{I}^{-i}}\big(u^{k\ast}(t,x)-g_{ik}(t,x)\big)+\varepsilon\geq\max\limits_{k\in\mathcal{I}^{-i}}\big(u^{k}(t,x)-g_{ik}(t,x)\big)+\varepsilon,
\end{equation}
since $u^{k\ast}\geq u^{k}$, and
\begin{equation}\label{weak4.5}
     -\left(\partial_t+\mathcal{L}\right)\phi(t,x)-f_i(t,x,(u^{k*})_{k=1,\dots,q}(t,x),\left(\sigma^{\mathtt{T}}D_x\right)\phi(t,x))\geq\varepsilon.
\end{equation}
By definition of $u^{i\ast}$, there exists a sequence $(t_m,x_m)_{m\geq 0}$ in $B_{\eta_\varepsilon}(t_0,x_0)$, such that $(t_m,x_m)\rightarrow (t_0,x_0)$ and $u^{i}(t_m,x_m)\rightarrow u^{i\ast}(t_0,x_0)$. Now let us fix $m$ and take the associated state process $X^{t_m,x_m}$ defined in (\ref{eds}) and define the stopping time $\theta^{m}$ as
\begin{equation}\label{stoppingtime}
     \theta^{m}:=(t_0+\eta_\varepsilon)\wedge\inf\left\{s\geq t_m: \left|X^{t_m,x_m}_{s}-x_0\right|\geq\eta_\varepsilon\right\}.
\end{equation}
Applying It\^o's lemma to $\phi(s,X^{t_m,x_m}_{s})$, it can be seen that
\begin{equation}
     \left(\phi(s,X^{t_m,x_m}_{s}),\left(\sigma^{\mathtt{T}}D_x\right)\phi(s,X^{t_m,x_m}_{s});\hspace{0.2cm} t_m\leq s\leq\theta^{m}\right)\notag
\end{equation}
is the solution of the BSDE with coefficient $-\left(\partial_s+\mathcal{L}\right)\phi(s,x)$, terminal time $\theta^{m}$ and terminal value $\phi(\theta^{m},X^{t_m,x_m}_{\theta^{m}})$. The idea is to compare this BSDE with the solution $(\underbar Y^{i,t_m,x_m}_{s})_{t_m\leq s\leq\theta^m}$ of the RBSDE with coefficient $f_i$, barrier $\max_{k\in\mathcal{I}^{-i}}\{u^{k}-g_{ik}\}$ and terminal condition $u^{i\ast}(\theta^m,X^{t_m,x_m}_{\theta^m})$. Note that by definition of $\theta^{m}$ and inequality (\ref{weak4.5}), we have
\begin{equation}
\begin{array}{ll}
-\left(\partial_s+\mathcal{L}\right)\phi\big(s,X^{t_m,x_m}_{s}\big)\geq f_i\left(s,X^{t_m,x_m}_{s},(u^{k*})_{k=1,\dots,q}(s,X^{t_m,x_m}_{s}),\left(\sigma^{\mathtt{T}}D_x\right)\phi(s,X^{t_m,x_m}_{s})\right)+\varepsilon & \\ \hspace{4.2cm}\geq f_i\left(s,X^{t_m,x_m}_{s},(u^{k})_{k=1,\dots,q}(s,X^{t_m,x_m}_{s}),\left(\sigma^{\mathtt{T}}D_x\right)\phi(s,X^{t_m,x_m}_{s})\right)+\varepsilon \notag
\end{array}
\end{equation}
for each $t_m\leq s\leq\theta^{m}$, where to reach the last inequality we use that $u^{\ast}\geq u$, the monotonicity property \ref{H2}-\ref{H2iv} and the Remark \ref{fimonotone}. It remains to compare the solution $\phi(s,X^{t_m,x_m}_{s})$ of the BSDE with the barrier $\max_{k\in\mathcal{I}^{-i}}\{u^{k}(s,X^{t_m,x_m}_{s})-g_{ik}(s,X^{t_m,x_m}_{s})\}\mathbf{1}_{[s<\theta^{m}]}+u^{i\ast}(s,X^{t_m,x_m}_{s})\mathbf{1}_{[s=\theta^{m}]}$ of the RBSDE for $t_m\leq s\leq\theta^m$. From inequality (\ref{weak33}) and definition of $\theta^{m}$ we derive that
\begin{equation}\label{barrier1}
\begin{array}{ll}
\phi(s,X^{t_m,x_m}_{s})\geq\max\limits_{k\in\mathcal{I}^{-i}}\left(u^{k\ast}(s,X^{t_m,x_m}_{s})-g_{ik}(s,X^{t_m,x_m}_{s})\right)+\varepsilon & \\ \hspace{2.2cm}\geq\max\limits_{k\in\mathcal{I}^{-i}}\left(u^{k}(s,X^{t_m,x_m}_{s})-g_{ik}(s,X^{t_m,x_m}_{s})\right)+\varepsilon \hspace{0.5cm} \text{for}\hspace{0.3cm} t_m\leq s<\theta^{m}.
\end{array}
\end{equation}
On the other hand, to show that the inequality holds at $\theta^m$, we recall that the minimum $(t_0,x_0)$ is strict and hence there exists $\gamma_{\varepsilon}$ such that
\begin{equation}
     \phi(t,x)-u^{i\ast}(t,x)\geq\gamma_\varepsilon \hspace{0.3cm}\text{on}\hspace{0.1cm} [0,T]\times\mathbb{R}^{r} \setminus B_{\eta_\varepsilon}({t_0,x_0}).\notag
\end{equation}
In particular, we have
\begin{equation}\label{barrier2}
     \phi(\theta^{m},X^{t_m,x_m}_{\theta^{m}})\geq u^{i\ast}(\theta^{m},X^{t_m,x_m}_{\theta^{m}})+\gamma_\varepsilon.
\end{equation}
Therefore, from (\ref{barrier1}), (\ref{barrier2}) and letting $\delta_\varepsilon:=\min(\varepsilon,\gamma_\varepsilon)$, we get
\begin{equation}
     \phi(s,X^{t_m,x_m}_{s})\geq\max\limits_{k\in\mathcal{I}^{-i}}\left(u^{k}(s,X^{t_m,x_m}_{s})-g_{ik}(s,X^{t_m,x_m}_{s})+\delta_{\varepsilon}\right)\mathbf{1}_{[s<\theta^{m}]}+\left(u^{i\ast}(\theta^{m},X^{t_m,x_m}_{\theta^{m}})+\delta_{\varepsilon}\right)\mathbf{1}_{[s=\theta^{m}]}\notag
\end{equation}
for $t_m\leq s\leq \theta^{m} \hspace{0.1cm}\text{a.s.}$. Thus, by the comparison result in Lemma \ref{compareBR}, we have
\begin{equation}
     \phi(s,X^{t_m,x_m}_{s})\geq \underbar Y^{i,t_m,x_m}_{s}+\delta_\varepsilon K \hspace{1cm}\text{for}\hspace{0.3cm} t_m\leq s\leq \theta^{m}\notag
\end{equation}
where $K$ is a positive constant which only depends on $T$ and the Lipschitz constant of $f_i$. In particular, for $t=t_m$, we have
\begin{equation}
\phi(t_m,x_m)\geq \underbar Y^{i,t_m,x_m}_{t_m}+\delta_\varepsilon K. \notag
\end{equation}
 Now, since $u^{i}(t_m,x_m)\rightarrow u^{i\ast}(t_0,x_0)$ and $\phi$ is continuous with $\phi(t_0,x_0)=u^{i\ast}(t_0,x_0)$, for $m$ sufficiently large we have both
\begin{equation}
|u^i(t_m,x_m)-u^{i\ast}(t_0,x_0)|\leq\frac{1}{4}\delta_\varepsilon K
\end{equation}
and
\begin{equation}
|u^{i\ast}(t_0,x_0)-\phi(t_m,x_m)|\leq\frac{\delta_\varepsilon K}{4},
\end{equation}
whence $|\phi(t_m,x_m)-u^{i}(t_m,x_m)|\leq \frac{1}{2}\delta_\varepsilon K$, and hence
\begin{equation}\label{fatality}
u^{i}(t_m,x_m)\geq \underbar Y^{i,t_m,x_m}_{t_m}+\frac{1}{4}\delta_\varepsilon K.
\end{equation}
But $u^{i\ast}(\theta^m,X^{t_m,x_m}_{\theta^m})\geq u^{i}(\theta^m,X^{t_m,x_m}_{\theta^m})$, then by comparison theorem $u^{i}(s,X^{t_m,x_m}_{s})=Y^{i,t_m,x_m}_s\leq \underbar Y^{i,t_m,x_m}_s$ for $t_m\leq s\leq\theta^m$. Thus, for $s=t_m$, we get $u^{i}(t_m,x_m)\leq \underbar Y^{i,t_m,x_m}_{t_m}$ that produces a contradiction with (\ref{fatality}). Therefore \eqref{subsolutionproperty} holds true and then also the viscosity subsolution property in $[0,T)\times \mathbb{R}^r$.
\ms

\noindent \underline{Step 2}: Viscosity super-solution property on $[0,T)\times \mathbb{R}^r$. 
\ms

\noindent Let $(t_0,x_0)\in[0,T)\times\mathbb{R}^{r}$ and $\phi\in C^{1,2}([0,T]\times\mathbb{R}^{r})$ be such that $\phi(t_0,x_0)=u_{\ast}^{i}(t_0,x_0)$ and $\phi(t,x)\leq u^{i}_{\ast}(t,x)$, for all $(t,x)\in[0,T]\times\mathbb{R}^{r}$. As stated above, we can suppose that the maximum is strict in $(t_0,x_0)$. Since by construction $u^{i}\geq\max_{k\in\mathcal{I}^{-i}}\left(u^{k}-g_{ik}\right)$, then it is easy to see that $u^{i}_{\ast}(t_0,x_0)\geq(\max_{k\in\mathcal{I}^{-i}} (u^{k}_{\ast}-g_{ik}))_{\ast}(t_0,x_0)$. Now, we show that
\begin{equation}
     -\left(\partial_t+\mathcal{L}\right)\phi(t_0,x_0)-f_i\big(t_0,x_0,(u^{k}_*)_{k=1,\dots,q}(t_0,x_0),\left(\sigma^{\mathtt{T}}D_x\right)\phi(t_0,x_0)\big)\geq 0.\notag
\end{equation}
Similar to the subsolution case, we shall proceed by contradiction, namely, suppose that 
\begin{equation}
     -\left(\partial_t+\mathcal{L}\right)\phi(t_0,x_0)-f_i\big(t_0,x_0,(u^{k}_*)_{k=1,\dots,q}(t_0,x_0),\left(\sigma^{\mathtt{T}}D_x\right)\phi(t_0,x_0)\big)<0,\notag
\end{equation}
then by Remark \ref{remark55} there exists $\varepsilon>0$ and $\eta_\varepsilon>0$ such that, for all $(t,x)\in B_{\eta_{\varepsilon}}(t_0,x_0)$, we have
\begin{equation}\label{weak4.8}
     -\left(\partial_t+\mathcal{L}\right)\phi(t,x)-f_i\big(t,x,(u^{k}_*)_{k=1,\dots,q}(t,x),\left(\sigma^{\mathtt{T}}D_x\right)\phi(t,x)\big)\leq-\varepsilon.
\end{equation}
Let $(t_m,x_m)_{m\geq 1}$ be a sequence in $B_{\eta_\varepsilon}(t_0,x_0)$ such that $(t_m,x_m)\rightarrow(t_0,x_0)$ and $u^{i}(t_m,x_m)\rightarrow u_{\ast}^{i}(t_0,x_0)$. We introduce the state process $X^{t_m,x_m}$ and define the stopping time $\theta^{m}$ as in (\ref{stoppingtime}). Next, we apply It\^o's formula to $\phi(s,X^{t_m,x_m}_{s})$ in order to obtain
\begin{equation}
\left(\phi(s,X^{t_m,x_m}_{s}),\left(\sigma^{\mathtt{T}}D_x\right)\phi(s,X^{t_m,x_m}_{s}); t_m\leq s\leq \theta^{m}\right)\notag
\end{equation}
is the solution of the BSDE associated with terminal time $\theta^{m}$, terminal value $\phi(\theta^{m},X^{t_m,x_m}_{\theta^{m}})$ and driver $( -\left(\partial_t+\mathcal{L}\right)\phi(s,X^{t_m,x_m}_{s}))_{s\in [t_m,\theta_m]}$. Then by definition of $\theta^{m}$ and inequality (\ref{weak4.8}), we get
\begin{equation}\label{weak4.9}
\begin{array}{ll}
   \hspace{-0.4cm} -\left(\partial_t+\mathcal{L}\right)\phi\left(s,X^{t_m,x_m}_{s}\right)\leq f_i\big(s,X^{t_m,x_m}_{s},(u^{k}_*)_{k=1,\dots,q}(s,X^{t_m,x_m}_{s}),\left(\sigma^{\mathtt{T}}D_x\right)\phi(s,X^{t_m,x_m}_{s})\big) -\varepsilon & \\ \hspace{3.8cm}\leq f_i\big(s,X^{t_m,x_m}_{s},(u^{k})_{k=1,\dots,q}(s,X^{t_m,x_m}_{s}),\left(\sigma^{\mathtt{T}}D_x\right)\phi(s,X^{t_m,x_m}_{s})\big)-\varepsilon
    \end{array}
\end{equation}
for $t_m\leq s\leq \theta^{m}$ a.s., where to reach the last inequality we use the monotonicity property  (H2)-(iv) and Remark \ref{fimonotone} and that $u^{j}\geq u^{j}_{\ast}$ for $j=1,\ldots,q$. It remains to compare the terminal conditions of the BSDEs with coefficients $ -\left(\partial_t+\mathcal{L}\right)\phi$ and $f_i$ respectively. Since the maximum $(t_0,x_0)$ is strict, there exists $\gamma_\varepsilon$ (which depends
on $\eta_\varepsilon$) such that $u_{\ast}^{i}(t,x)\geq\phi(t,x)+\gamma_\varepsilon$ on $[0,T]\times\mathbb{R}^{r}\setminus B_{\eta_\varepsilon}(t_0,x_0)$, which implies
\begin{equation}
\phi(\theta_m,X^{t_m,x_m}_{\theta^{m}} )\leq u^{i}_{\ast}(\theta^{m},X^{t_m,x_m}_{\theta^{m}})-\gamma_\varepsilon. \notag
\end{equation}
Thus using inequality (\ref{weak4.9}) and the comparison result for BSDEs, Lemma \ref{lemma46}, we derive that
\begin{equation}
     \phi(s,X^{t_m,x_m}_{s})\leq \bar Y^{i,t_m,x_m}_{s},\notag \hspace{0.3cm} \text{for}\hspace{0.2cm}t_m\leq s\leq\theta^m
\end{equation}
and therefore, in $s=t_m$, we have $\phi(t_m,x_m)\leq \bar Y^{i,t_m,x_m}_{t_m}$. As above mentioned, we can assume that $m$ is sufficient large so that $|\phi(t_m,x_m)-u^{i}(t_m,x_m)|\leq\frac{\delta_{\varepsilon}K}{2}$. We thus get
\begin{equation}
  u^{i}(t_m,x_m)-\frac{\gamma_\varepsilon K}{2}\leq \phi(t_m,x_m)\leq \bar Y^{i,t_m,x_m}_{t_m}\notag
\end{equation}
and hence
\begin{equation}\label{weak4.10}
u^i(t_m,x_m)<\bar Y^{i,t_m,x_m}_{t_m}.
\end{equation}
But $u^{i}_{\ast}(\theta^m,X^{t_m,x_m}_{\theta^m})\leq u^{i}(\theta^m,X^{t_m,x_m}_{\theta^m})$, then by Lemma \ref{lemma46} we get $\bar Y^{i,t_m,x_m}_{s}\leq Y^{i,t_m,x_m}_{s}=u^{i}(s,X^{t_m,x_m}_s)$ for $t_m\leq s\leq\theta^m$, and thus $\bar Y^{i,t_m,x_m}_{t_m}\leq u^{i}(t_m,x_m)$, which is a contradiction with (\ref{weak4.10}). Therefore the viscosity supersolution property in $[0,T)\times \mathbb{R}^r$ holds true.
\medskip

\noindent \underline{Step 3}: Subsolution property at $(T,x)$.
\medskip

\noindent We now show that for any $i=1,...,m$, 
\begin{equation}
\min\big\{u^{i\ast}(T,x_0)-h_{i}(x_0);\hspace{0.2cm} u^{i\ast}(T,x_0)-\max\limits_{j\in\mathcal{I}^{-i}}\left(u^{j\ast}-g_{ij})\right)^{\ast}(T,x_0)\big\}\leq 0.\notag
\end{equation}
We  follow here the same idea as in Bouchard \cite{Bouchard2009Stochastic} (see also Theorem 1 in Hamad\`ene and Morlais  \cite{hamadene2013viscosity}). We reason by contradiction, namely, we assume that
\begin{equation}\label{contradiction}
\min\big\{u^{i\ast}(T,x_0)-h_{i}(x_0);\hspace{0.2cm} u^{i\ast}(T,x_0)-\max\limits_{j\in\mathcal{I}^{-i}}\left(u^{j\ast}-g_{ij}\right)^{\ast}(T,x_0)\big\}=2\varepsilon>0.
\end{equation}
Let $(t_k,x_k)$ be a sequence in $[0,T)\times\mathbb{R}^{k}$ such that
\begin{equation}\label{approxima}
(t_k,x_k)\rightarrow (T,x_0) \hspace{0.3cm}\text{and}\hspace{0.3cm} u^{i}(t_k,x_k)\rightarrow u^{i\ast}(T,x_0) \hspace{0.3cm} \text{as} \hspace{0.2cm} k\rightarrow \infty.
\end{equation}
Since $u^{i\ast}$ is \textit{u.s.c} and of polynomial growth, we can find a sequence $\left(\varphi_{n}\right)_{n\geq 0}$  of functions of $C^{1,2}([0,T]\times\mathbb{R}^{k})$ and neighborhood $B_n$ of $(T,x_0)$ such that $\varphi^{n}\rightarrow u^{i\ast}$, and hence from the inequality (\ref{contradiction}) we have
\begin{equation}\label{T1}
\min\big\{\varphi^{n}(t,x)-h_i(x);\hspace{0.2cm} \varphi^{n}(t,x)-\max\limits_{j\in\mathcal{I}^{-i}}\left(u^{j\ast}-g_{ij}\right)^{\ast}(t,x)\big\}\geq \varepsilon \hspace{0.3cm} \text{for all}\hspace{0.2cm} (t,x)\in B_n,
\end{equation}
for $n$ large enough. On the other hand, after possibly passing to a sub-sequence of $(t_k,x_k)_{k\geq 1}$ we can assume that the previous inequality holds on $B^{k}_{n}:=[t_k,T]\times B(x_k,\delta^{k}_{n})$ for some $\delta^{k}_{n}\in (0,1)$ small enough in such a way that $B^{k}_{n}\subset B_n$. Since $u^{i\ast}$ is locally bounded (recall it has polynomial growth), there exists $\zeta>0$ such that $\left|u^{i\ast}\right|\leq\zeta$ on $B_n$. We can then assume that $\varphi^{n}\geq-2\zeta$ on $B_n$. Next we define
\begin{equation}
\tilde{\varphi}^{n}_{k}(t,x):=\varphi^{n}(t,x)+\frac{4\zeta\left|x-x_k\right|^{2}}{(\delta^{k}_{n})^{2}}+\sqrt{T-t}.\notag
\end{equation}
Note that $\tilde{\varphi}^{n}_{k}\geq\varphi^{n}$ and
\begin{equation}\label{T2}
\left(u^{i\ast}-\tilde{\varphi}^{n}_{k}\right)(t,x)\leq-\zeta \hspace{0.3cm} \text{for} \hspace{0.2cm} (t,x)\in[t_k,T]\times\partial B(x_k,\delta^{k}_{n}).
\end{equation}
Since $\partial_t(\sqrt{T-t})\rightarrow -\infty$ as $t\rightarrow T$, we can choose $t_k$ close enough to $T$ to ensure that
\begin{equation}\label{T3}
-\left(\partial_t+
\mathcal{L}\right)\tilde{\varphi}^{n}_{k}(t,x)\geq 0 \hspace{0.3cm}\text{on}\hspace{0.3cm} B^{k}_{n}.
\end{equation}
Next let us consider the following stopping times
\begin{equation}\label{stoppingtime1}
\theta^{k}_{n}:=\inf\{s\geq t_k: (s,X^{t_k,x_k}_{s})\in B^{k^c}_{n}\}\wedge T \hspace{0.3cm}
\end{equation}

and
\begin{equation}\label{stoppingtime2}
\vartheta^{\varepsilon}_k:=\inf\{s\geq t_k, u^{i}(s,X^{t_k,x_k}_{s})\leq\max_{j\in\mathcal{I}^{-i}}(u^{j}(s,X^{t_k,x_k}_{s})-g_{ij}(s,X^{t_k,x_k}_{s}))+\frac{\varepsilon}{4}\}\wedge T
\end{equation}
where $B^{k^c}_{n}$ is the complement of $B^{k}_{n}$. 

First note that for a subsequence $\{k\}$, $\P[\vartheta^{\varepsilon}_k>t_k]=1$. Actually from 
\eqref{contradiction}, we have 
\begin{align*}
u^{i\ast}(T,x_0)& \ge\max\limits_{j\in\mathcal{I}^{-i}}\left(u^{j\ast}-g_{ij}\right)^{\ast}(T,x_0)+2\varepsilon\\
&\ge \max\limits_{j\in\mathcal{I}^{-i}}\left(u^{j}-g_{ij}\right)^{\ast}(T,x_0)+2\varepsilon .
\end{align*}
Therefore taking into account of \eqref{approxima}, at least for a subsequence, for any $k\ge 1$,
\begin{align*}
u^{i\ast}(t_k,x_k) \ge \max\limits_{j\in\mathcal{I}^{-i}}\left(u^{j}-g_{ij}\right)(t_k,x_k)+\varepsilon.
\end{align*}
Now let us stick to this subsequence. If $\P[\vartheta^{\varepsilon}_k=t_k]>0$, then by the c\`adl\`ag property of the processes which define $\vartheta^{\varepsilon}_k$ we have 
$u^{i}({t_k,x_k})\leq\max_{j\in\mathcal{I}^{-i}}(u^{j}({t_k,x_k})-g_{ij}({t_k,x_k}))+\frac{\varepsilon}{4}$, which contradicts the previous inequality and then the claim is valid. 

On the other hand the property which characterizes the jumps of $Y^i$ in the definition \eqref{RBSDEM5.4}, implies that on $[t_k, \vartheta^{\varepsilon}_k]$ the process $Y^i$ is continuous and $dK^i_s=0$ for $s\in [t_k, \vartheta^{\varepsilon}_k]$. Applying now It\^o's formula to the process $(\tilde{\varphi}^{n}_{k}(s,X_s))_{s\in [t_k, \theta^{k}_{n}\wedge \vartheta^{\varepsilon}_k]}$ and taking expectation, we obtain

\begin{equation}\label{T6}
\begin{array}{ll}
\hspace{-0.7cm}\tilde{\varphi}^{n}_{k}(t_k,x_k)=\mathbb{E}\bigg[\tilde{\varphi}^{n}_{k}(\theta^{k}_{n}\wedge\vartheta^{\varepsilon}_k,X^{t_k,x_k}_{\theta^{k}_{n}\wedge\vartheta^{\varepsilon}_k})-\displaystyle\int^{\theta^{k}_{n}\wedge\vartheta^{\varepsilon}_k}_{t_k}\left(\partial_t+\mathcal{L}\right)\tilde{\varphi}^{n}_{k}(s,X^{t_k,x_k}_{s})ds\bigg]  & \\ & \\ \hspace{1.1cm} \geq\mathbb{E}\big[\tilde{\varphi}^{n}_{k}(\theta^{k}_{n},X^{t_k,x_k}_{\theta^{k}_{n}})\mathbf{1}_{[\theta^{k}_{n}\leq\vartheta^{\varepsilon}_k]}+\tilde{\varphi}^{n}_{k}(\vartheta^{\varepsilon}_k,X^{t_k,x_k}_{\vartheta^{\varepsilon}_k})\mathbf{1}_{[\theta^{k}_{n}>\vartheta^{\varepsilon}_k]}\big] \hspace{1cm}\text{by (\ref{T3})} & \\ & \\ \hspace{1.1cm}=\mathbb{E}\big[\big\{\tilde{\varphi}^{n}_{k}(\theta^{k}_{n},X^{t_k,x_k}_{\theta^{k}_{n}})\mathbf{1}_{[\theta^{k}_{n}<T]}+\tilde{\varphi}^{n}_{k}(\theta^{k}_{n},X^{t_k,x_k}_{\theta^{k}_{n}})\mathbf{1}_{[\theta^{k}_{n}=T]}\big\}\mathbf{1}_{[\theta^{k}_{n}\leq\vartheta^{\varepsilon}_k]}+\tilde{\varphi}^{n}_{k}(\vartheta^{\varepsilon}_k,X^{t_k,x_k}_{\vartheta^{\varepsilon}_k})\mathbf{1}_{[\theta^{k}_{n}>\vartheta^{\varepsilon}_k]}\big]& \\ & \\ \hspace{1.1cm}\geq\mathbb{E}\Big[\big\{\big(u^{i\ast}(\theta^{k}_{n},X^{t_k,x_k}_{\theta^{k}_{n}})+\zeta\big)\mathbf{1}_{[\theta^{k}_{n}<T]}+\big(h_i(T,X^{t_k,x_k}_{T})+\varepsilon\big)\mathbf{1}_{[\theta^{k}_{n}=T]}\big\}\mathbf{1}_{[\theta^{k}_{n}\leq\vartheta^{\varepsilon}_k]}  & \\ & \\ \hspace{2.7cm}+\Big\{\max\limits_{j\in\mathcal{I}^{-i}}\big(u^{j\ast}-g_{ij}\big)^{\ast}(\vartheta^{\varepsilon}_k,X^{t_k,x_k}_{\vartheta^{\varepsilon}_k})+\varepsilon\Big\}\mathbf{1}_{[\theta^{k}_{n}>\vartheta^{\varepsilon}_k]}\Big] \hspace{1cm}\text{by (\ref{T2}) and (\ref{T1})}  & \\

 & \\ \hspace{1.1cm}\geq\mathbb{E}\Big[\big\{\big(u^{i}(\theta^{k}_{n},X^{t_k,x_k}_{\theta^{k}_{n}})+\zeta\big)\mathbf{1}_{[\theta^{k}_{n}<T]}+\big(h_i(T,X^{t_k,x_k}_{T})+\varepsilon\big)\mathbf{1}_{[\theta^{k}_{n}=T]}\big\}\mathbf{1}_{[\theta^{k}_{n}\leq\vartheta^{\varepsilon}_k]} & \\ 
 
 & \\ \hspace{2.7cm} +\Big\{\max\limits_{j\in\mathcal{I}^{-i}}\big(u^{j}(\vartheta^{\varepsilon}_k,X^{t_k,x_k}_{\vartheta^{\varepsilon}_k})-g_{ij}(\vartheta^{\varepsilon}_k,X^{t_k,x_k}_{\vartheta^{\varepsilon}_k})\big)+\varepsilon\Big\}\mathbf{1}_{[\theta^{k}_{n}>\vartheta^{\varepsilon}_k]}\Big] \notag
 \end{array}
 \end{equation}

 \begin{equation}
 \begin{array}{ll}
 \hspace{1.1cm}\geq\mathbb{E}\Big[\big\{\big(u^{i}(\theta^{k}_{n},X^{t_k,x_k}_{\theta^{k}_{n}})+\zeta\big)\mathbf{1}_{[\theta^{k}_{n}<T]}+\big(h_i(T,X^{t_k,x_k}_{T})+\varepsilon\big)\mathbf{1}_{[\theta^{k}_{n}=T]}\big\}\mathbf{1}_{[\theta^{k}_{n}\leq\vartheta^{\varepsilon}_k]} +\big\{u^{i}(\vartheta^{\varepsilon}_k,X^{t_k,x_k}_{\vartheta^{\varepsilon}_k})+\frac{3\varepsilon}{4}\big\}\mathbf{1}_{[\theta^{k}_{n}>\vartheta^{\varepsilon}_k]}\Big]\hspace{1cm}\text{by (\ref{stoppingtime2})}  & \\ & \\ \hspace{1.1cm}\geq\mathbb{E}\big[{u^{i}(\theta^{k}_{n}\wedge\vartheta^{\varepsilon}_k,X^{t_k,x_k}_{\theta^{k}_{n}\wedge\vartheta^{\varepsilon}_k})}\big]+\left(\zeta\wedge\frac{3\varepsilon}{4}\right) & \\ & \\ \hspace{1.1cm}=\mathbb{E}\big[u^{i}(t_k,x_k)\big]-\mathbb{E}\bigg[\displaystyle\int^{\theta^{k}_{n}\wedge\vartheta^{\varepsilon}_k}_{t_k}f_{i}(s,X^{t_k,x_k}_{s},(u^{k})_{k=1,\dots,q}(s,X^{t_k,x_k}_{s}),Z^{i,t_k,x_k}_{s})ds\bigg]+\left(\zeta\wedge\frac{3\varepsilon}{4}\right)
\end{array}
\end{equation}
where the last equality is due to the fact that the process $Y^{i}_{.}=u^{i}(\cdot,X_{.})$, stopped at time $\theta^{k}_{n}\wedge\vartheta_{k}$, solves a RBSDE system of the type (\ref{RBSDEM5.4}) with  data given by $((f_i)_{i\in\mathcal{I}}, (h_i)_{i\in\mathcal{I}},(g_{ij} )_{i\in\mathcal{I}})$, and the last inequality is obtained by monotonicity property of $f_i$ and since $u^{j\ast}\geq u^{j}$ for $j\in\mathcal{I}^{-i}$. Besides, note that by definition of $\theta^{k}_{n}\wedge\vartheta_{k}$ we have $dK^{i,t,x}=0$ on $[t_k,\vartheta^{\varepsilon}_k]$. Next, we have that both $(u^{j})_{j=1,\ldots,m}$ and $(t,x)\rightarrow\big\|Z^{i,t,x}_{\cdot}\big\|_{\mathcal{H}^{2,d}}(t,x)$ are of polynomial growth. Thus by Assumption (H2)-(i),(iii) and inequality (\ref{5.33}) we deduce that
\begin{equation}\label{T4}
\lim\limits_{k\rightarrow\infty}\mathbb{E}\bigg[\displaystyle\int^{\theta^{k}_{n}\wedge\vartheta^{\varepsilon}_k}_{t_k}f_{i}(s,X^{t_k,x_k}_{s},(u^{k})_{k=1,\dots,q}(s,X^{t_k,x_k}_{s}),Z^{i,t_k,x_k}_{s})ds\bigg]=0,
\end{equation}
and hence taking the limit in both hand sides of the inequality (\ref{T6}) as $k\rightarrow\infty$ yields
\begin{equation}
\begin{array}{ll}
\varphi^{n}(T,x_0)=\lim\limits_{k\rightarrow\infty}\left[\varphi^{n}(t_k,x_k)+\sqrt{T-t_k}\right]=\lim\limits_{k\rightarrow\infty}\tilde{\varphi}^{n}_{k}(t_k,x_k) & \\ &\\ \hspace{3.1cm}\geq\lim\limits_{k\rightarrow\infty}u^{i}(t_k,x_k)+\left(\zeta\wedge\frac{3\varepsilon}{4}\right)=u^{i\ast}(T,x_0)+\left(\zeta\wedge\frac{3\varepsilon}{4}\right).
\end{array}
\end{equation}
Therefore, taking $n$ large enough and recalling that $\varphi^{n}\rightarrow u^{i\ast}$ pointwisely, we get a contradiction.  Thus for any $x\in\mathbb{R}^{k}$ and $i\in\mathcal{I}$ we have
\begin{equation}
\min\Big\{u^{i\ast}(T,x)-h_{i}(x);\hspace{0.2cm} u^{i\ast}(T,x)-\max\limits_{j\in\mathcal{I}^{-i}}\left(u^{j\ast}-g_{ij}\right)^{\ast}(T,x)\Big\}\leq 0.
\end{equation}
which is the claim.
\ms

\noindent \underline{Step 4}: Supersolution property at $(T,x_0)$.
\ms

\noindent We are going to show that 
\begin{equation}
\min\Big\{u^{i}_{\ast}(T,x_0)-h_{i}(x_0);\hspace{0.2cm} u^{i}_{\ast}(T,x_0)-\big(\max\limits_{j\in\mathcal{I}^{-i}}\left(u^{j}_{\ast}(T,x_0)-g_{ij}(T,x_0)\right)\big)_{\ast}\Big\}\geq 0.
\end{equation}
Let $(t_k,x_k)_{k\geq 1}$ be a sequence in $[0,T)\times\mathbb{R}^{d}$ such that
\begin{equation}\label{T5}
(t_k,x_k)\rightarrow (T,x_0)\hspace{0.2cm}\text{and}\hspace{0.2cm} u^{i}(t_k,x_k)\rightarrow u^{i}_{\ast}(T,x_0)\hspace{0.4cm} \text{as}\hspace{0.4cm} k\rightarrow\infty.
\end{equation}
Since $u^{i}(t,x)$ is deterministic, we have from the definition of $u^{i}$ that
\begin{equation}\label{T7}
\begin{array}{ll}
u^{i}(t_k,x_k)=\mathbb{E}\bigg[h_i(X^{t_k,x_k}_{T})+\displaystyle\int^{T}_{t_k}f_{i}(s,X^{t_k,x_k}_{s},(u^{k})_{k=1,\dots,q}(s,X^{t_k,x_k}_{s}),Z^{i,t_k,x_k}_{s})ds+K^{i}_{T}-K^{i}_{t_k}\bigg] & \\ & \\ \hspace{1.7cm}\geq \mathbb{E}\bigg[h_i(X^{t_k,x_k}_{T})+\displaystyle\int^{T}_{t_k}f_{i}(s,X^{t_k,x_k}_{s},(u^{k})_{k=1,\dots,q}(s,X^{t_k,x_k}_{s}),Z^{i,t_k,x_k}_{s})ds\bigg]
\end{array}
\end{equation}
where we have used that $dK^{i,t,x}\geq 0$ on $[t_k,T]$. Next taking the limit in both hand sides as $k\rightarrow\infty$, using that $h_i$ is continuous and  arguing similarly to (\ref{T4}) we have
\begin{align}
u^{i}_{\ast}(T,x_0)&\geq\lim_{k\rightarrow\infty}\mathbb{E}\bigg[h_i(X^{t_k,x_k}_{T})+\int^{T}_{t_k}f_{i}(s,X^{t_k,x_k}_{s},(u^{k})_{k=1,\dots,q}(s,X^{t_k,x_k}_{s}),Z^{i,t_k,x_k}_{s})ds\bigg]\nonumber \\&=\mathbb{E}\big[h_i(X^{T,x_0}_{T})\big]=h_i(x_0),\nonumber 
\end{align}
that is, $u^{i}_{\ast}(T,x_0)\geq h_i(x_0)$. On the other hand, setting $\tau_k=(T+t_k)/2$, considering the RBSDE \eqref{RBSDEM5.4} on $[t_k,\tau_k]$, taking expectation to obtain 
\begin{equation}
\begin{array}{ll}
u^{i}(t_k,x_k)\geq\mathbb{E}\bigg[u^{i}(\tau_k,X^{t_k,x_k}_{\tau_k})+\displaystyle\int^{\tau_k}_{t_k}f_{i}(s,X^{t_k,x_k}_{s},(u^{k})_{k=1,\dots,q}(s,X^{t_k,x_k}_{s}),Z^{i,t_k,x_k}_{s})ds\bigg]& \\ & \\ \hspace{1.6cm} \geq\mathbb{E}\bigg[\max\limits_{j\in\mathcal{I}^{-i}}\big(u^{j}(\tau_k,X^{t_k,x_k}_{\tau_k})-g_{ij}(\tau_k,X^{t_k,x_k}_{\tau_k})\big)+\int^{\tau_k}_{t_k}f_{i}(s,X^{t_k,x_k}_{s},(u^{k})_{k=1,\dots,q}(s,X^{t_k,x_k}_{s}),Z^{i,t_k,x_k}_{s})ds\bigg]
& \\ & \\ \hspace{1.6cm} \geq\mathbb{E}\bigg[(\max\limits_{j\in\mathcal{I}^{-i}}\big(u^{j}_*-g_{ij})\big)_*(\tau_k,X^{t_k,x_k}_{\tau_k})+\int^{\tau_k}_{t_k}f_{i}(s,X^{t_k,x_k}_{s},(u^{k})_{k=1,\dots,q}(s,X^{t_k,x_k}_{s}),Z^{i,t_k,x_k}_{s})ds\bigg]
\end{array}
\end{equation}
since $dK^{i,t,x}\geq 0$ and $u^{i}(\tau_k,X^{t_k,x_k}_{\tau_k})\ge \max\limits_{j\in\mathcal{I}^{-i}}\big(u^{j}(\tau_k,X^{t_k,x_k}_{\tau_k})-g_{ij}(\tau_k,X^{t_k,x_k}_{\tau_k})\big)$. It implies that 
\begin{align}
\lim_{k\rightarrow \infty}u^{i}(t_k,x_k)&\geq\liminf_{k\rightarrow \infty}\mathbb{E}\bigg[(\max\limits_{j\in\mathcal{I}^{-i}}\big(u^{j}_*-g_{ij})\big)_*(\tau_k,X^{t_k,x_k}_{\tau_k})\nonumber\\
&\hspace{2cm}+\int^{\tau_k}_{t_k}f_{i}(s,X^{t_k,x_k}_{s},(u^{k})_{k=1,\dots,q}(s,X^{t_k,x_k}_{s}),Z^{i,t_k,x_k}_{s})ds\bigg]\nonumber\\
&
\geq\mathbb{E}\bigg[\liminf_{k}(\max\limits_{j\in\mathcal{I}^{-i}}\big(u^{j}_*-g_{ij})\big)_*(\tau_k,X^{t_k,x_k}_{\tau_k})\nonumber\\
&\hspace{2cm}+\int^{\tau_k}_{t_k}f_{i}(s,X^{t_k,x_k}_{s},(u^{k})_{k=1,\dots,q}(s,X^{t_k,x_k}_{s}),Z^{i,t_k,x_k}_{s})ds\bigg]\nonumber\\
&\ge 
(\max\limits_{j\in\mathcal{I}^{-i}}\big(u^{j}_*-g_{ij})\big)_*(t_k,x_k).\nonumber
\end{align}
The second inequality stems from Fatou's Lemma while the third one is due to the fact that 
$(\max\limits_{j\in\mathcal{I}^{-i}}\big(u^{j}_*-g_{ij})\big)_*$ is lower semicontinuous and by \eqref{T5}, at least for a subsequence, 
$((\tau_k,X^{t_k,x_k}_{\tau_k}))_k\rightarrow (T,x)$ $\P-a.s.$. Thus 
\begin{equation}
\min\Big\{u^{i}_{\ast}(T,x_0)-h_{i}(x_0);\hspace{0.2cm} u^{i}_{\ast}(T,x_0)-\big(\max\limits_{j\in\mathcal{I}^{-i}}\left(u^{j}_{\ast}(T,x_0)-g_{ij}(T,x_0)\right)\big)_{\ast}\Big\}\geq 0
\end{equation}
which is the claim. The proof is now complete. $\qed$
\begin{REM} If the switching costs $g_{ij}$ are continuous, conditions \eqref{defsursolviscoenT} and \eqref{defsoussolviscoenT}, read respectively as:
\begin{equation}\notag
\min\big\{v^{i\ast}(T,x_0)-h_{i}(x_0);\hspace{0.2cm} u^{i\ast}(T,x_0)-\max\limits_{j\in\mathcal{I}^{-i}}\left(v^{j\ast}-g_{ij}\right)(T,x_0)\big\}\leq 0
\end{equation}
and 
\begin{equation}\notag
\min\Big\{v^{i}_{\ast}(T,x_0)-h_{i}(x_0);\hspace{0.2cm} v^{i}_{\ast}(T,x_0)-\big(\max\limits_{j\in\mathcal{I}^{-i}}\left(v^{j}_{\ast}-g_{ij}\right)\big)(T,x_0)\Big\}\geq 0. 
\end{equation}Therefore $v^{i}_{\ast}(T,x_0)\ge h_{i}(x_0)$ and by the non free-loop property one 
deduces that $v^{i\ast}(T,x_0)\le h_{i}(x_0)$ which implies that $v^{i}(T,x_0)=h_{i}(x_0)$. For more details one can see e.g. Hamad\`ene and Morlais  \cite{hamadene2013viscosity}. $\qed$
\end{REM}

\end{document}